\documentclass[12pt]{article}
\usepackage{fancyhdr,amsmath, psfrag, graphicx, amsfonts, layout, color,enumerate}
\usepackage[numbers,sort]{natbib}
\let\cite=\citet

\usepackage{ifpdf}
\ifpdf 
\DeclareGraphicsExtensions{.pdf,.pdftex_t,.png,.jpg}
\else 
\DeclareGraphicsExtensions{.eps,.ps,.pst}
\fi
\usepackage[english]{babel}

\font\fontemail=cmtt12
\font\fontauthors=cmcsc10 scaled \magstep1
 at 8truept
 at 8truept

\newtheorem{Th}{Theorem}
\newtheorem{Def}[Th]{Definition}
\newtheorem{Cor}[Th]{Corollary}
\newtheorem{Ex}[Th]{Example}
\newtheorem{Lem}[Th]{Lemma}
\newtheorem{Prop}[Th]{Proposition}
\newtheorem{Rem}[Th]{Remark}

\newcommand{\pff}{\noindent {\sc Proof.}\ }

\def\QED{\hfill\vrule height 1.5ex width 1.4ex depth -.1ex \vskip 10pt}

\newtheorem{prerem}[Th]{Remark}
\newenvironment{rem}{\begin{prerem}\rm }{\end{prerem}}

\def\cl{\mathop{\rm cl}\nolimits}

\def\ie{\emph{i.e.~}}

\def\intr{\mathop{\rm int}\nolimits}

\newcommand{\lpref}{\smash{\raisebox{3.5pt}{\!\!\!\begin{tabular}{c}$\hskip-4pt\scriptstyle\longleftarrow$ \\[-7pt]{\rm pref}\end{tabular}\!\!}}}
\newcommand{\petitlpref}{\smash{\raisebox{2.5pt}{\!\!\!\begin{tabular}{c}$\hskip-2pt\scriptscriptstyle\longleftarrow$ \\[-9pt]{$\scriptstyle\hskip 1pt\rm pref$}\end{tabular}\!\!}}}

\def\ordPart#1{\mathop{\rm \sum _{#1}\uparrow}\nolimits}
\def\Proba{\mathop{\bf {~\!\! P}}\nolimits}

\def\stirling2 #1#2{\left\{\begin{matrix} #1\\#2\end{matrix}\right\}}
\def\suff{\mathop{\rm suff}\nolimits}

\newcommand\1{\leavevmode\hbox{\rm \small1\kern-0.35em\normalsize1}}
\newcommand\ind[1]{\1_{\{#1\}}}

\newcommand{\R}{\rm I\!R}
  
\newcommand{\N}{{{\rm I}\!{\rm N}}}
\newcommand{\Q}{{\rm Q}\kern-.65em{}^{{}_/}}

\def\1{{\bf 1}}

\def\g#1{\mathbb #1}
\def\rond#1{\mathcal #1}
\def\R{\g R}

\def\N{\g N}
\def\Z{\g Z}

\begin{document}

\begin{center}
\LARGE{\bf {Variable length Markov chains\\ and dynamical sources}}
\end{center}

\begin{picture}(400,320)
\put(0,0){\includegraphics[height=300pt]{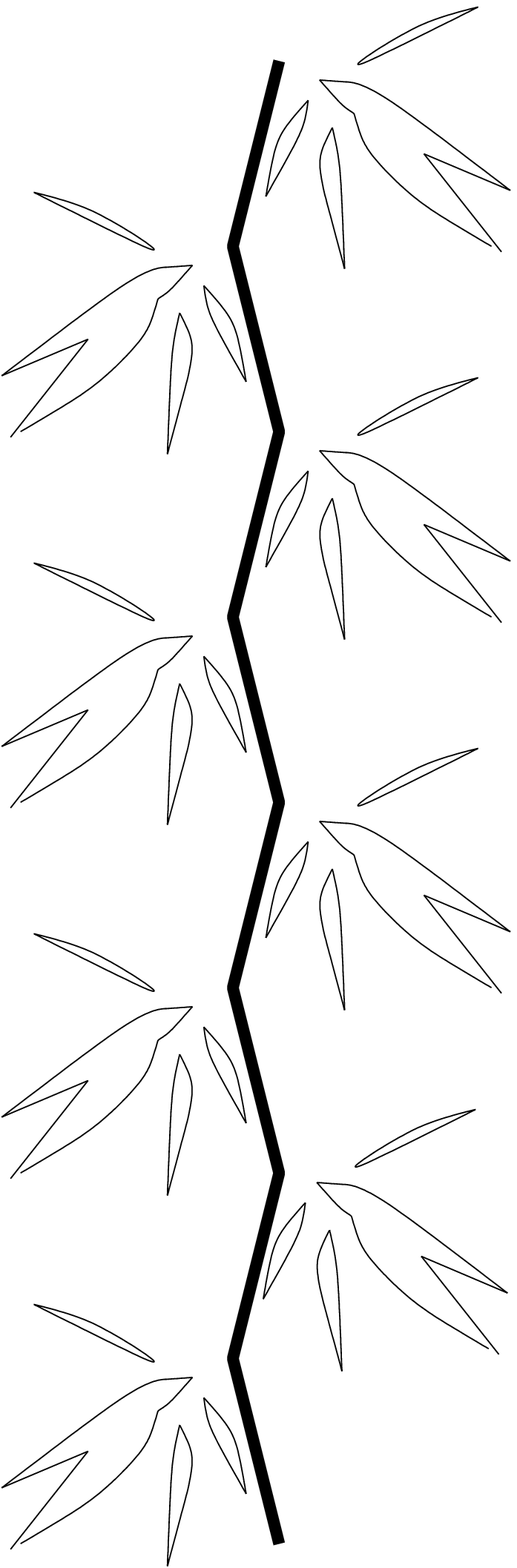}}
\put(150,150){\begin{minipage}{40truemm}\begin{center}
{\fontauthors
Peggy C\'enac

\medskip
Brigitte Chauvin

\medskip
Fr\'ed\'eric Paccaut

\medskip
Nicolas Pouyanne}

\vskip 20pt
{\it 17 juillet 2010}
\end{center}\end{minipage}}
\put(315,-40){\includegraphics[height=370pt,width=110pt]{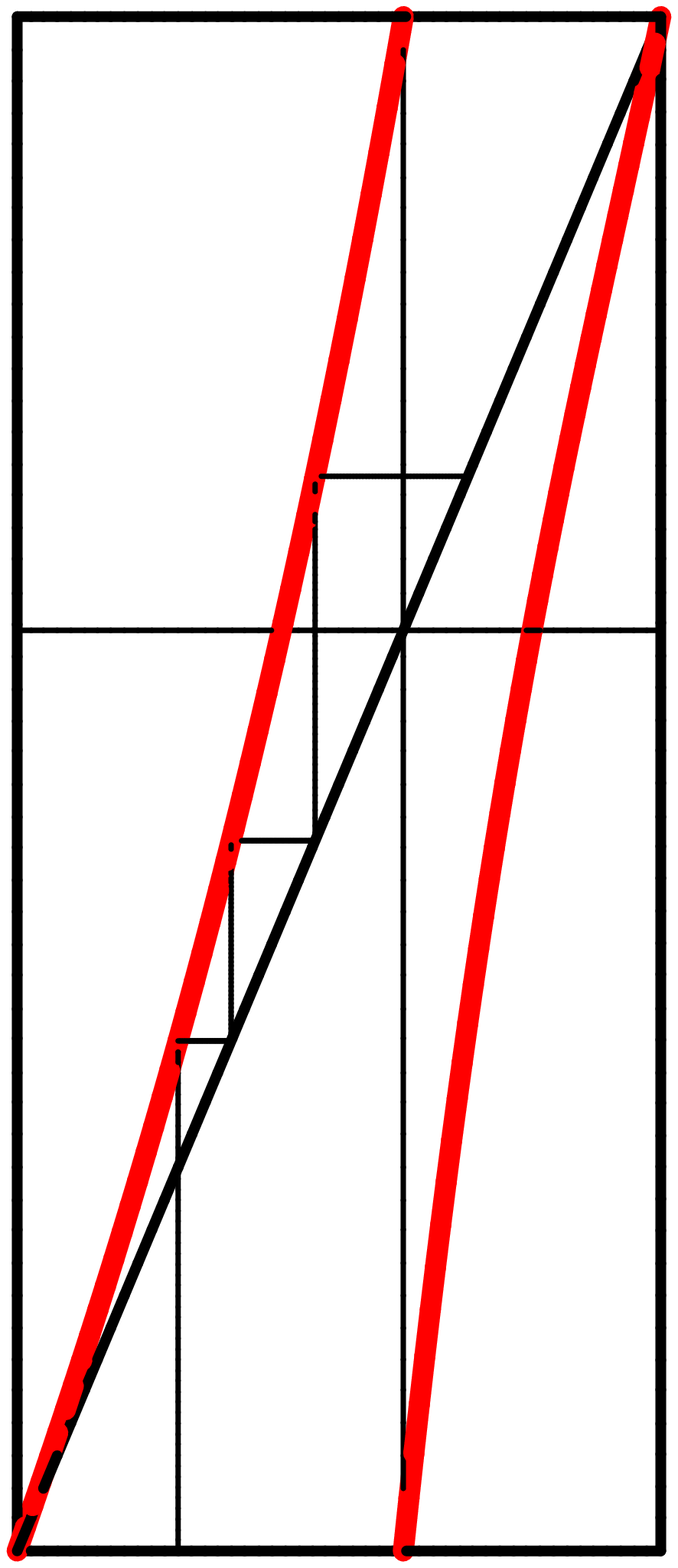}}
\end{picture}

\vskip 20pt
{\small
{\bf Abstract}

\medskip

Infinite random sequences of letters can be viewed as stochastic chains or as strings produced by a source, in the sense of information theory. 
The relationship between Variable Length Markov Chains (VLMC) and probabilistic dynamical sources is studied. We establish a probabilistic frame for context trees and VLMC and we prove that any VLMC is a dynamical source for which we explicitly build the mapping. On two examples, the ``comb'' and the ``bamboo blossom'', we find a necessary and sufficient condition for the existence and the unicity of a stationary probability measure for the VLMC. These two examples are detailed in order to provide the associated Dirichlet series as well as the generating functions of word occurrences. 


}

\tableofcontents

\section{Introduction}
\label{sec:Intro}

Our objects of interest are infinite random sequences of letters.
One can imagine DNA sequences (the letters are $A$, $C$, $G$, $T$), bits sequences (the letters are $0$, $1$) or any  random sequence on a finite alphabet. 
Such a sequence can be viewed as a stochastic chain or as a string produced by a source, in the sense of information theory. 
We study this relation for the so-called Variable Length Markov Chains (VLMC).

\medskip

From now on, we are given a finite alphabet $\rond A$.
An infinite random sequence of letters is often considered as a chain $(X_n)_{n\in\Z}$, \ie an
$\rond A^\g Z$-valued random variable.
The $X_n$ are the \emph{letters} of the chain.
Equivalently such a chain can be viewed as a random process $(U_n)_{n\in \g N}$ that takes values in the set
$\rond L:=\rond A^{-\N}$ of left-infinite words\footnote{In the whole text, $\N$ denotes the set of nonnegative integers.} and that grows by addition of a letter on the right at each step of discrete time. The $\rond L$-valued processes we consider are Markovian ones.
The evolution from $U_n = \dots X_{-1}X_0X_1 \dots X_n$ to $U_{n+1} = U_nX_{n+1}$ is described by
the transition probabilities $\Proba(U_{n+1} = U_n\alpha | U_n)$, $\alpha\in \rond A$.

\medskip

In the context of chains, the point of view has mainly been a statistical one until now, 
going back to \cite{Harris} who speaks of \emph{chains of infinite order} to express the fact that the production of a new letter depends on a finite but unbounded number of previous letters. \cite{Comets} and \cite{Gallo-Garcia} deal with chains of infinite memory. \cite{Rissanen}  introduces a class of models where the transition from the word $U_n$ to the word $U_{n+1} = U_nX_{n+1}$ depends on $U_n$ through a finite suffix of $U_n$ and he calls this relevant part of the past a \emph{context}. Contexts can be stored as the leaves of a so-called \emph{context tree} so that the model is entirely defined by a family of probability distributions indexed by the leaves of a context tree. In this paper, Rissanen develops a near optimal universal data compression algorithm for long strings generated by non independent information sources. The name VLMC is due to \cite{Buhlmann}. It emphasizes the fact that  the length of memory needed to predict the next letter is a not necessarily bounded  function of the sequence $U_n$. An overview on VLMC can be found in \cite{Galves-Locherbach}.

\medskip

We give in Section~\ref{sec:Model} a complete probabilistic definition of VLMC.
Let us  present here a foretaste, relying on the particular form of the transition probabilities
$\Proba(U_{n+1} = U_n\alpha | U_n)$.
Let $\rond T$ be a \emph{saturated tree} on $\rond A$, which means that every internal node of the tree -- \ie a word on $\rond A$ -- has exactly $|\rond A|$ children. With each leaf $c$ of the tree, also called a \emph{context}, is associated a probability distribution $q_c$ on $\rond A$.
The basic fact is that any left-infinite sequence can thus be ``plugged in'' a unique context of the tree $\rond T$: any $U_n$ can be uniquely written $U_n = \dots \overline{c}$, where, for any word
$c=\alpha_1\cdots\alpha _N$, $\overline{c}$ denotes the reversed word
$\overline{c}=\alpha_N\cdots\alpha _1$.
In other terms, for any $n$, there is a unique context $c$ in the tree $\rond T$ such that $\overline{c}$ is a suffix of $U_n$; this word is denoted by $c=\lpref(U_n)$.
We define the VLMC associated with these data as the $\rond L$-valued homogeneous Markov process whose transition probabilities are, for any letter $\alpha\in\rond A$,
$$
\Proba(U_{n+1} = U_n\alpha | U_n) = q_{\petitlpref(U_n)}(\alpha).
$$
When the tree is finite, the final letter process $(X_n)_{n\geq 0}$ is an ordinary Markov chain whose order is the height of the tree.
The case of infinite trees is more interesting, providing concrete examples of non Markov chains.
\begin{figure}[!h]
\begin{picture}(400,210)
\put(30,35){\includegraphics[height=160pt]{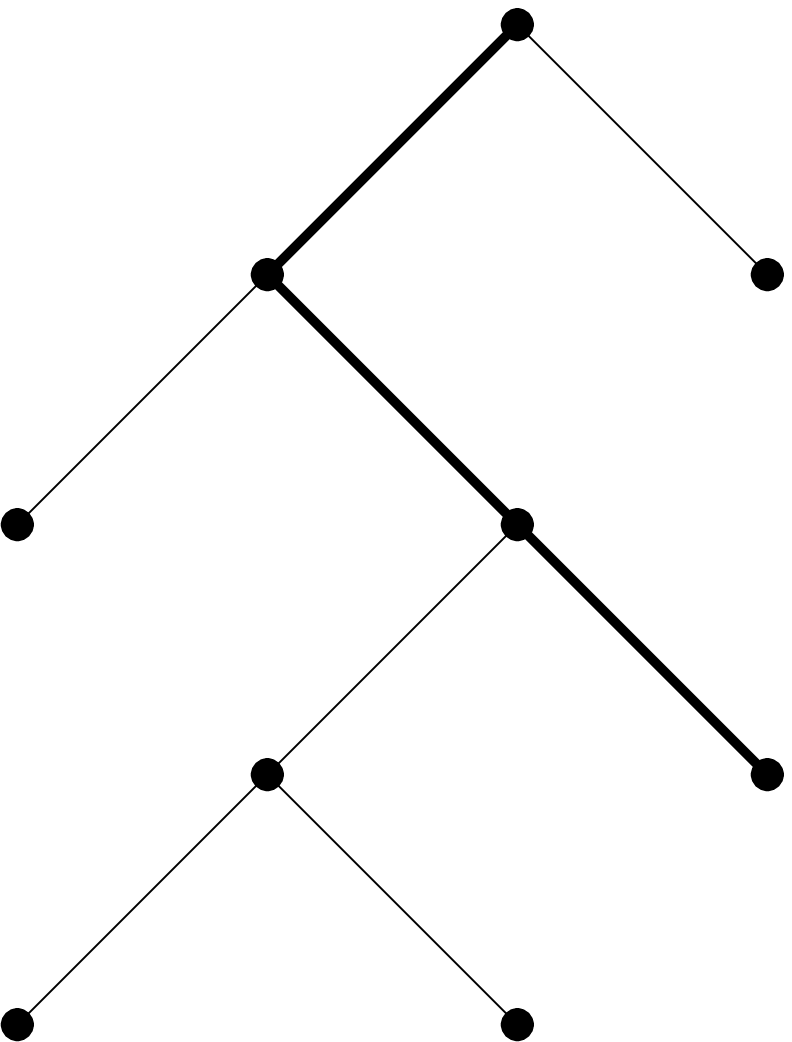}}
\put(198.8,13.5){\includegraphics[height=203.3pt]{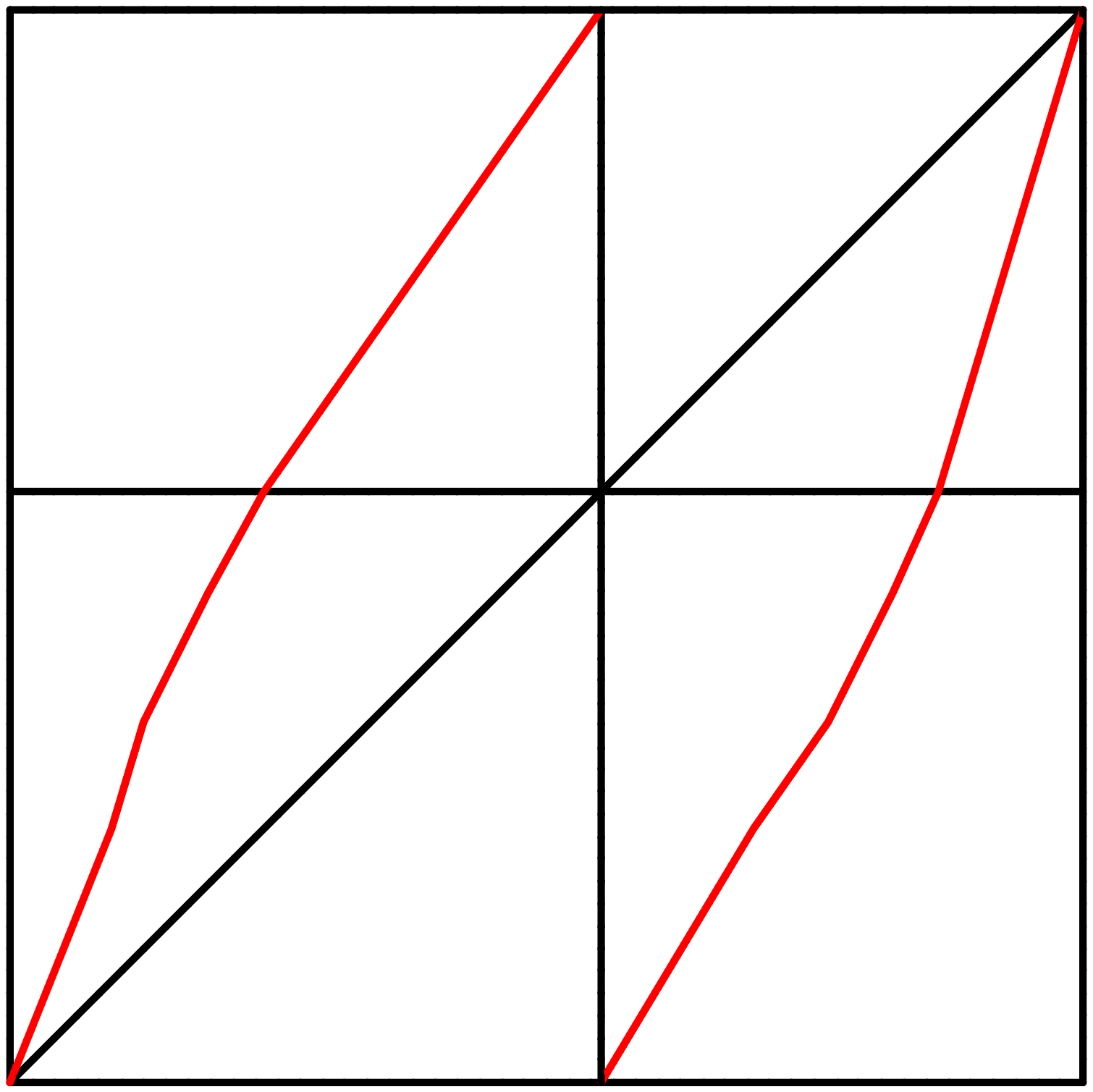}}
\put(22,21){$q_{0100}$}
\put(100,21){$q_{0101}$}
\put(140,60){$q_{011}$}
\put(24,99){$q_{00}$}
\put(144,139){$q_{1}$}

\put(280,32){\line(0,1){120.5}}
\put(280,152.5){\line(1,0){58}}
\put(338.3,153){\line(0,-1){121}}
\put(338.1,84.3){\line(-1,0){68}}
\put(270.2,32){\line(0,1){53}}

\put(280,20){$x$}
\put(332,20){$Tx$}
\put(256,20){$T^2x$}

\put(222,11){\line(1,0){156.8}}
\put(222.2,8){\line(0,1){6}}
\put(308.4,8){\line(0,1){6}}
\put(378.7,8){\line(0,1){6}}

\put(260,0){$I_0$}
\put(340,0){$I_1$}
\end{picture}
\caption{\label{exempleIntro}example of probabilized context tree (on the left) and its corresponding dynamical system
(on the right).}
\end{figure}

In the example of Figure~\ref{exempleIntro}, the context tree is finite of height $4$ and, for instance,
$\Proba(U_{n+1} = U_n0 | U_n=\cdots 0101110) = q_{011}(0)$ because $\lpref (\cdots 0101{\bf{110}})={\bf{011}}$ (read the word $\cdots 0101110$ right-to-left and stop when finding a context).

\medskip

In information theory, one considers that words are produced by a \emph{probabilistic source} as developed in Vall\'ee and her group papers  (see \cite{ClementFlajoletVallee} for an overview).
In particular, a \emph{probabilistic dynamical source} is defined by a coding function $\rho : [0,1] \to\rond A$,
a mapping $T: [0,1] \to [0,1]$ having suitable properties and a probability measure $\mu$ on $[0,1]$.
These data being given, the dynamical source produces the $\rond A$-valued random process 
$(Y_n)_{n\in\g N}:= (\rho(T^n\xi))_{n\in\g N}$, where $\xi$ is a $\mu$-distributed random variable on $[0,1]$.
On the right side of Figure~\ref{exempleIntro}, one can see the graph of some $T$, a subdivision of $[0,1]$ in two subintervals
$I_0=\rho ^{-1}(0)$ and $I_1=\rho ^{-1}(1)$ and the first three real numbers $x$, $Tx$ and $T^2x$, where $x$ is a realization of the random variable $\xi$. The right-infinite word corresponding to this example has $010$ as a prefix.

We prove in Theorem \ref{equivalenceProcessus} that every stationary VLMC is a dynamical source.
More precisely, given a stationary VLMC,  $(U_n)_{n\in \g N}$ say, we construct explicitly a dynamical source $(Y_n)_{n \in \g N}$ such
that the letter processes $(X_n)_{n\in\g N}$ and $(Y_n)_{n\in\g N}$ are symmetrically distributed, which means that for any finite word $w$ of length $N+1$, $\Proba(X_0\dots X_N=w) = \Proba(Y_0\dots Y_N=\overline{w})$.
In Figure~\ref{exempleIntro}, the dynamical system together with Lebesgue measure on $[0,1]$ define a probabilistic source that corresponds to the stationary VLMC defined by the drawn probabilized context tree.

\medskip

The previous result is possible only when the VLMC is stationary. The question of existence and unicity of a stationary distribution arises naturally. We give a complete answer in two particular cases (Proposition \ref{mesureStationnairePeigne} and Proposition \ref{mesureStationnaireBambou} in Section \ref{sec:ex}) and we propose some tracks for the general case. Our two examples are called the ``infinite comb'' and the ``bamboo blossom''; they can be visualized in Figures \ref{figPeigneInfini} and \ref{figBambou}, respectively page \pageref{figPeigneInfini} and page \pageref{figBambou}. 
Both have an infinite branch so that the letter process of the VLMC is non Markovian. They provide quite concrete cases of infinite order chains where the study can be completely handled. We first exhibit a necessary and sufficient condition for existence and unicity of a stationary measure. Then the dynamical system is explicitly built and drawn. In particular, for some suitable data values, one gets in this way examples of intermittent sources. 

Quantifying and visualizing repetitions of patterns is another natural question arising in combinatorics on words. \emph{Tries}, \emph{suffix tries} and \emph{digital search trees} are usual convenient tools. The analysis of such structures relies on the generating functions of the word occurrences and on the Dirichlet series attached to the sources. In both examples, these computations are performed.

\medskip

The paper is organized as follows. Section~\ref{sec:Model} is devoted to the precise definition of variable length Markov chains. In Section~\ref{sec:stat_vlmc} the main result Theorem~\ref{equivalenceProcessus} is established.
In Section~\ref{sec:ex}, we complete the paper with our two detailed examples: ``infinite comb'' and ``bamboo blossom''. The last section gathers some prospects and open problems.

\section{Context trees and variable length Markov chains}
\label{sec:Model}

In this section, we first define \emph{probabilized context trees}; then we associate with a probabilized context tree a so-called \emph{variable length Markov chain (VLMC)}.

\subsection{Words and context trees}\label{trees}

Let $\rond A$ be a finite alphabet, \ie a finite ordered set. Its cardinality is denoted by $|\rond A|$.
For the sake of shortness, our results in the paper are given for the alphabet $\rond A = \{ 0,1\}$ but they remain true for any finite alphabet. Let
$$\rond W=\displaystyle\bigcup_{n\geq 0}\rond A^n$$ 
be the set of all finite words over~$\rond A$. The concatenation of two words $v=v_1\dots v_M$ and $w= w_1\dots w_N$  is $vw = v_1\dots v_Mw_1\dots w_N$. The empty word is denoted by  $\emptyset$. 
Let 
$$\rond L=\rond A^{-\N}$$
be the set of  left-infinite sequences over $\rond A$ and  
$$\rond R=\rond A^{\N}$$
be the set of  right-infinite sequences over $\rond A$.  If $k$ is a nonnegative integer and if $w=\alpha  _{-k}\cdots \alpha _0$ is any finite word on $\rond A$, the reversed word is denoted by 
\[\overline w = \alpha_0\cdots\alpha_{-k}.\]
The \emph{cylinder based on $w$} is defined as the set of all left-infinite sequences having $w$ as a suffix:
$$\rond L w=\{ s\in\rond L,~\forall j\in\{ -k,\cdots ,0\},~s_j=\alpha _j\}.$$ 
By extension, the reversed sequence of $s = \cdots\alpha_{-1}\alpha_0\in\rond L$ is $\overline s = \alpha_0\alpha_{-1}\cdots \in\rond R$. The set $\rond L$ is equipped with the $\sigma$-algebra generated by all cylinders based on finite words. The set $\rond R$ is equipped with the $\sigma$-algebra generated by all cylinders $w\rond R= \{ r\in\rond R, w$ is a prefix of $r \}$.

\medskip

Let $\rond T$ be a tree, \emph{i.e.\@} a subset of  $\rond W$ satisfying two conditions: 
\begin{itemize}
\item[$\bullet$] $\emptyset \in \rond T$
\item[$\bullet$]  $\forall u,v \in \rond W$,  $uv\in \rond T$  $\Longrightarrow$ $ u\in \rond T$.
\end{itemize}
This corresponds to the definition of rooted planar trees in algorithmics. Let $\rond C^F(\rond T)$ be the set of \emph{finite leaves} of $\rond T$, \ie the nodes of $\rond T$ without any descendant:
$$
\rond C^F(\rond T) = \{ u\in \rond T, \forall j\in\rond A, uj\notin \rond T\}.
$$
An infinite word $u\in \rond R$ such that any finite prefix of $u$ belongs to $\rond T$ is called an \emph{infinite leaf} of $\rond T$. Let us denote the set of infinite leaves of $\rond T$ by
$$
\rond C^I(\rond T)=\{u\in\rond R, \forall v \hbox{ prefix of } u, v\in \rond T\}.
$$
Let $\rond C(\rond T)=\rond C^F(\rond T)\cup\rond C^I(\rond T)$ be the set of all \emph{leaves} of  $\rond T$. The set $\rond T\setminus \rond C^F(\rond T)$ is constituted by the \emph{internal nodes} of $\rond T$. When there is no ambiguity, $\rond T$ is omitted and we simply write $\rond C, \rond C^F$ and $\rond C^I$.
\begin{Def}
A tree is \emph{saturated} when each internal node $w$ has exactly $|\rond A|$ children, namely the set  $\{ w\alpha, \alpha \subset \rond A\} \subset \rond T$.
\end{Def}

\begin{Def} {\bf (Context tree)}

A \emph{context tree} is a saturated tree having a finite or denumerable set of leaves. The leaves are called \emph{contexts}.
\end{Def}
\begin{Def} {\bf (Probabilized context tree)}

A \emph{probabilized context tree} is a pair
$$\left( \rond T, (q_c)_{c\in\rond C(\rond T) }\right)$$
where $\rond T$ is a context tree over $\rond A$ and $(q_c)_{c\in\rond C(\rond T)}$ is a family of probability measures on $\rond A$, indexed by the denumerable set $\rond C(\rond T)$ of all leaves of $\rond T$.\
\end{Def}
\begin{Ex}
See  Figure~\ref{exempleIntro} for an example of finite probabilized context tree with five contexts.  See Figure~\ref{figPeigneInfini} for an example of infinite probabilized context tree, called the infinite comb.
\end{Ex}

\begin{Def}
A subset $\rond K$ of $\rond W \cup \rond R$ is a \emph{cutset} of the complete $|\rond A|$-ary tree when both following conditions hold

(i) no word of $\rond K$ is a prefix of another word of $\rond K$

(ii) $\forall r\in \rond R$, $\exists  u\in \rond K$, $u$ prefix of $r$.
\end{Def}
Condition (i) entails unicity in (ii). Obviously a tree $\rond T$ is saturated if and only if the set of its leaves $\rond C$ is a cutset.  Take a saturated tree, then
\begin{equation}\label{rcutset}
\forall r\in\rond R, \hbox{ either } r\in\rond C^I\hbox{ or } \exists ! u \in\rond W,~u\in\rond C^F,~u\hbox{ prefix of } r.
\end{equation}
This can also be said on left-infinite sequences:
\begin{equation}\label{lcutset}
\forall s\in\rond L, \hbox{ either } \overline s\in\rond C^I \hbox{ or } \exists ! w\in\rond W,~\overline w\in\rond C^F,~w\hbox{ suffix of } s.
\end{equation}
In other words:
\begin{equation}\label{partition}
\rond L = \bigcup_{\overline s\in\rond C^I}\{ s\} \cup \bigcup_{\overline w\in\rond C^F} \rond L w.
\end{equation}
This partition of $\rond L$ will be extensively used in the sequel. Both cutset properties (\ref{rcutset}) and  (\ref{lcutset}) will be used in the paper, on $\rond R$ for trees, on $\rond L$ for chains. Both orders of reading will be needed. 
\begin{Def} {\bf (Prefix function)}
Let $\rond T$ be a saturated tree and $\rond C$ its set of contexts.
For any $s \in \rond L$, $\lpref(s)$ denotes the unique context $\alpha_1\dots\alpha_N$ such that $s = \dots \alpha_N\dots\alpha_1$. The map 
\[\lpref: \rond L \to \rond C \]
is called the \emph{prefix function}. For technical reasons, this function is extended to 
\[\lpref: \rond L\cup \rond W \to \rond T \]
in the following way:
\begin{itemize}
\item if $\overline{w}\in \rond T$ then $\lpref(w)=\overline{w}$;
\item if $\overline{w}\in \rond W \setminus \rond T$ then $\lpref(w)$ is the unique context $\alpha_1\dots\alpha_N$ such that $w$ has $\alpha_N\dots\alpha_1$ as a suffix.
\end{itemize}
\end{Def}
Note that the second item of the definition is also valid when $\overline{w}\in \rond C$. Moreover $\lpref(w)$ is always a context except when $\overline{w}$ is an internal node.

\subsection{VLMC associated with a context tree}
\label{VLMC}

\begin{Def}\label{defVLMC} {\bf (VLMC)}

Let $\left( \rond T, (q_c)_{c\in\rond C}\right)$ be a probabilized context tree. The associated \emph{Variable Length Markov Chain} (VLMC) is the order $1$ Markov chain $(U_n)_{n\geq 0}$ with state space $\rond L$, defined by the transition probabilities
\begin{equation}\label{transition}
\forall n\geq 0,~\forall \alpha\in\rond A,~\Proba\left( U_{n+1}=U_n\alpha |U_n\right)=q_{\petitlpref (U_n)}\left(\alpha\right).
\end{equation}
\end{Def}
\begin{Rem}
As usually, we speak of \emph{the} Markov chain defined by the transition probabilities (\ref{transition}), because these data together with the distribution of $U_0$ define a unique $\rond L$-valued random process (see for example \cite{Revuz}).
\end{Rem}
The rightmost letter of the sequence $U_n\in\rond L$ will be denoted by $X_n$ so that
$$\forall n\geq 0,~U_{n+1}=U_nX_{n+1}.$$
The final letter process $(X_n)_{n\geq 0}$ is \emph{not} Markov as soon as the context tree has at least one infinite context. As already mentioned in the introduction, when the tree is finite, $(X_n)_{n\geq 0}$ is a Markov chain whose order is the height of the tree, \ie the length of its longest branch. The vocable VLMC is somehow confusing but commonly used.

\begin{Def} {\bf (SVLMC)}
\label{def:svlmc}
Let $(U_n)_{n\geq 0}$ be a VLMC. When a stationary probability measure on $\rond L$ exists and when it is the initial distribution, we  say that $(U_n)_{n\geq 0}$ is a \emph{Stationary Variable Length Markov Chain} (SVLMC).
\end{Def}
\begin{Rem}
In the literature, the name VLMC is usually applied to the chain $(X_n)_{n\in \Z}$. There exists a natural bijective correspondence between $\rond A$-valued chains $(X_n)_{n\in \Z}$ and $\rond L$-valued processes $(U_n = U_0X_1\dots X_n, n\geq 0)$. Consequently, finding a stationary probability for the chain $(X_n)_{n\in \Z}$ is equivalent to finding a stationary probability for the process $(U_n)_{n\geq 0}$.
\end{Rem}

\section{Stationary variable length Markov chains}
\label{sec:stat_vlmc}

The existence and the unicity of a stationary measure for two examples of VLMC will be established in Section~\ref{sec:ex}. In the present section, we assume that a stationary measure $\pi$ on $\rond L$ exists and we consider a $\pi$-distributed VLMC. In the preliminary Section~\ref{ssec:probaStationnaires},  we show how the stationary probability of finite words can be expressed as a function of the data and the values of $\pi$ at the tree nodes. In Section~\ref{ssec:systDynAssocie}, the main theorem is proved.

\subsection{General facts on stationary probability measures}
\label{ssec:probaStationnaires}

For the sake of shortness, when $\pi$ is a stationary probability for a VLMC, we write $\pi(w)$ instead of $\pi(\rond L w)$, for any $w\in\rond W$:
\begin{equation}
\label{notationpiw}
\pi(w)=\Proba(U_0 \in \rond L w)=\Proba(X_{-(|w|-1)}\ldots X_0=w).
\end{equation}

{\bf Extension of notation $q_u$ for internal nodes}.

The VLMC is defined by its context tree $\rond T$ together with a family $(q_c)_{c\in\rond C}$ of probability measures on $\rond A$ indexed by the contexts of the tree.
When $u$ is an internal node of the context tree, we extend the notation $q_{u}$ by
\begin{equation}
\label{probacondinternalnodes}
q_{u}(\alpha )=\left\{
\begin{array}{l}
\displaystyle\frac{\pi (\overline u\alpha )}{\pi (\overline u)}{\rm ~if~}\pi (\overline u)\neq 0\\ \\
0{\rm ~if~}\pi (\overline u)=0
\end{array}
\right.
\end{equation}
for any $\alpha\in\rond A$. Thus, in any case, $\pi$ being stationary, $\pi (\overline u\alpha )=\pi (\overline u)q_{u}(\alpha )$ as soon as $\overline u$ is an
internal node of the context tree.
With this notation, the stationary probability of any cylinder can be expressed by the following simple
Formula~(\ref{probaStationnaireCylindres}).

\begin{Lem}
\label{lemProbaStationnaireCylindres}
Consider a SVLMC defined by a probabilized context tree and let $\pi$ denote any
stationary probability measure on $\rond L$.
Then,

{\emph{(i)}}
for any finite word $w\in\rond W$ and for any letter $\alpha\in\rond A$,
\begin{equation}
\label{probaStationnaireCalcul}
\pi (w\alpha )=\pi (w)q_{\petitlpref (w)}(\alpha ).
\end{equation}

{\emph{(ii)}}
For any finite word $w=\alpha _1\dots\alpha _N\in\rond W$,
\begin{equation}
\label{probaStationnaireCylindres}
\pi (w)=\prod _{k=0}^{N-1}q_{\petitlpref (\alpha _1\dots \alpha _k)}(\alpha _{k+1})
\end{equation}
(if $k=0$, $\alpha _1\dots \alpha _k$ denotes the empty word $\emptyset$, $\lpref(\emptyset) = \emptyset$, $q_{\emptyset}(\alpha)=\pi(\alpha)$ and $\pi(\emptyset) = \pi(\rond L) = 1$).
\end{Lem}

\pff
{\emph{(i)}} If $\overline w$ is an internal node of the context tree, then $\lpref(w) = \overline{w}$ and the formula comes directly from the definition of $q_{\overline{w}}$.
If not, $\pi (w\alpha )=\pi (U_1\in\rond L {w\alpha})$ by stationarity;
because of Markov property,
\[\pi (w\alpha )=\Proba (U_0\in\rond L w)\Proba(U_1\in\rond L {w\alpha}|U_0\in\rond L w)
=\pi (w)q_{\petitlpref (w)}(\alpha ).\]
Finally, {\emph{(ii)}} follows from {\emph{(i)}} by a straightforward induction.
\QED

\begin{Rem}
\label{remPi10=Pi01}
When $\rond A=\{ 0,1\}$ and $\pi$ is any stationary probability of a SVLMC, then, for any
natural number $n$, $\pi (10^n)=\pi (0^n1)$.
Indeed, on one hand, by disjoint union, $\pi (0^n)=\pi (0^{n+1})+\pi (10^n)$.
On the other hand, by stationarity,
\begin{eqnarray*}
\pi (0^n)&=&\Proba (X_1\dots X_n=0^n)=\Proba (X_0\dots X_{n-1}=0^n)\\
&=&\Proba (X_0\dots X_n=0^{n+1})+\Proba (X_0\dots X_n=0^n1)=\pi (0^{n+1})+\pi (0^n1).
\end{eqnarray*}
These equalities lead to the result. Of course, symmetrically, $\pi (01^n)=\pi (1^n0)$ under the same assumptions.
\end{Rem}

\subsection{Dynamical system associated with a VLMC}
\label{ssec:systDynAssocie}
We begin with a general presentation of a probabilistic dynamical source in Section~\ref{sssec:generalProbaSource}. Then we build step by step partitions of the interval $[0,1]$ (Section~\ref{sssection:subdivisions}) and a mapping (Section~\ref{sectionDefinitionT}) based on the stationary measure of a given SVLMC. It appears in Section~\ref{sssec:proprietesT} that this particular mapping keeps Lebesgue measure invariant. All these arguments combine to provide in the last Section~\ref{sssection:VLMCisDS} the proof of Theorem \ref{equivalenceProcessus} which allows us to see a VLMC as a dynamical source.

\medskip

In the whole section, $I$ stands for the real interval $[0,1]$ and the Lebesgue measure of a Borelian  $J$ will be denoted by $|J|$. 
\subsubsection{General probabilistic dynamical sources}
\label{sssec:generalProbaSource}

\begin{figure}[!h]
\begin{picture}(300,220)
\thicklines
\put(105,20){\includegraphics[width=76truemm]{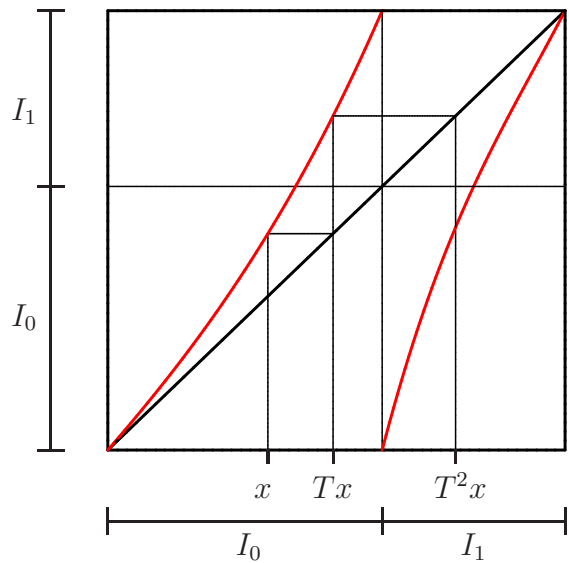}}
\put(126.6,18){\line(1,0){173}}
\put(175,5){$I_0$}
\put(260,5){$I_1$}
\put(126.6,13){\line(0,1){10}}
\put(230.4,13){\line(0,1){10}}
\put(299.6,13){\line(0,1){10}}

\put(105,44.5){\line(0,1){167}}
\put(90,92){$I_0$}
\put(90,170){$I_1$}
\put(100,44.8){\line(1,0){10}}
\put(100,144.7){\line(1,0){10}}
\put(100,211.3){\line(1,0){10}}

\put(187.2,40){\line(0,1){5}}
\put(211.9,40){\line(0,1){5}}
\put(258.1,40){\line(0,1){5}}
\put(182,27){$x$}
\put(204,27){$Tx$}
\put(250,27){$T^2x$}
\end{picture}
\caption{\label{4figorbite}the graph of a mapping $T$, the intervals $I_0$ and $I_1$ that code the interval $I$ by the alphabet $\rond A=\{ 0,1\}$ and the first three points of the orbit of an $x\in I$ by the corresponding dynamical system.}
\end{figure}

Let us present here the classical formalism of probabilistic dynamical sources (see \cite{ClementFlajoletVallee}). It is defined by four elements:
\begin{itemize}
\item a topological partition of $I$ by intervals $(I_\alpha)_{\alpha \in \rond A}$,
\item a coding function $\rho : I \to \rond A$, such that, for each letter $\alpha$, the restriction of $\rho$ to $I_{\alpha}$ is equal to $\alpha$,
\item a mapping $T: I\to I$, 
\item a probability measure $\mu$ on $I$.
\end{itemize}
Such a dynamical source defines an $\rond A$-valued random process $(Y_n)_{n\in\g N}$ as follows. Pick a random real number $x$ according to the measure $\mu$. The mapping $T$ yields the orbit $(x, T(x), T^2(x), \ldots)$ of $x$. Thanks to the coding function,  this defines the right-infinite  sequence $\rho(x) \rho( T(x)) \rho(T^2(x)) \cdots$  whose letters are $Y_n:=\rho(T^n(x))$ (see Figure~\ref{4figorbite}).

For any finite word $w=\alpha _0\dots\alpha _N\in\rond W$, let
$$B_w=\bigcap _{k=0}^NT^{-k}I_{\alpha _k}$$
be the Borelian set of real numbers $x$ such that the sequence
$(Y_n)_{n\in\g N}$
has $w$ as a prefix.
Consequently, the probability that the source emits a sequence of symbols starting with the pattern $w$ is equal to $\mu(B_w)$. When the initial probability measure $\mu$ on $I$ is $T$-invariant, the dynamical source generates a stationary
$\rond A$-valued random process which means that for any $n\in\g N$, the random variable $Y_n$ is
$\mu$-distributed.

\begin{figure}[h]
\begin{picture}(400,210)
\thicklines
\put(8,10){\includegraphics[height=200pt]{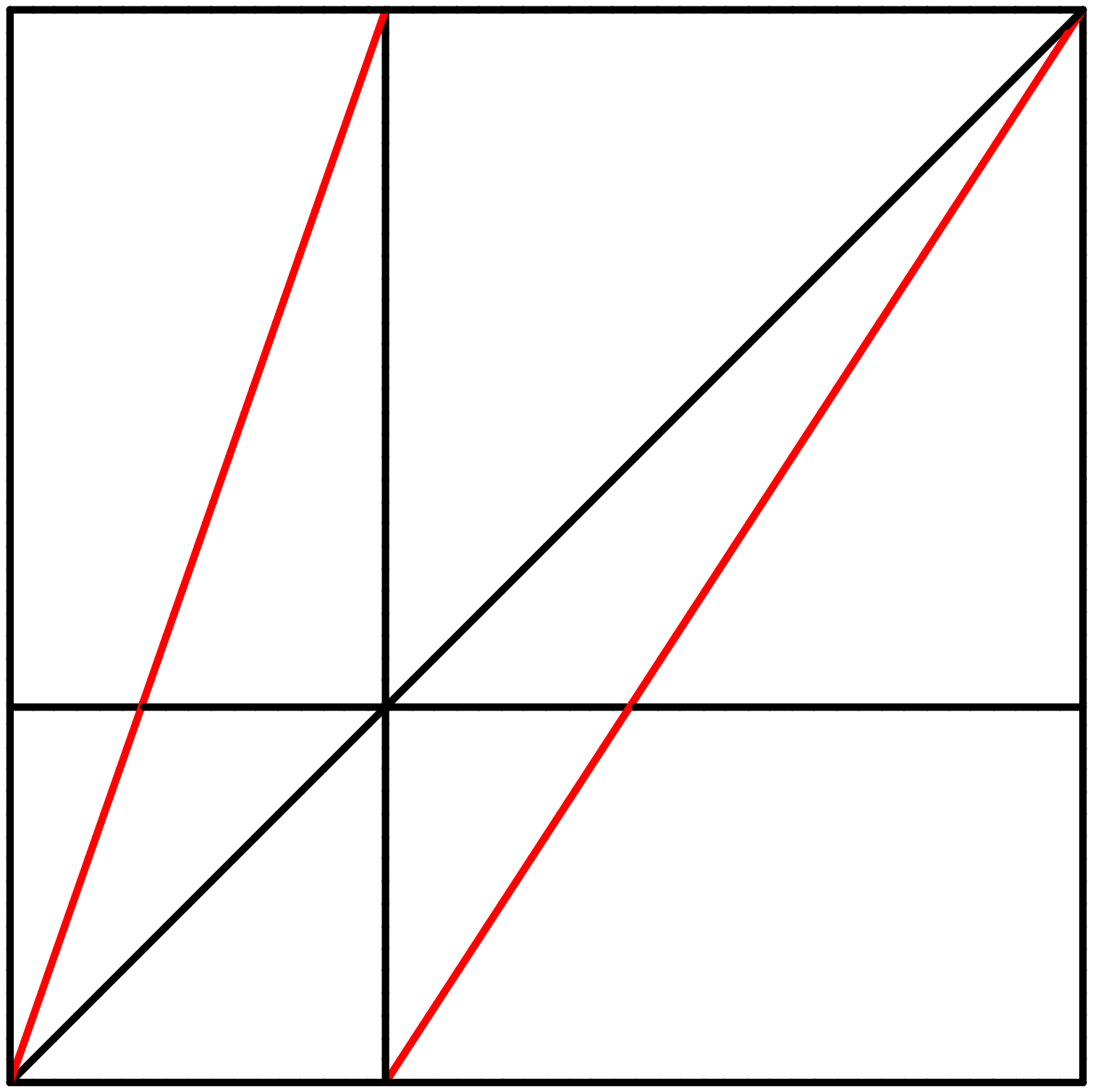}}
\put(200,10){\includegraphics[height=200pt]{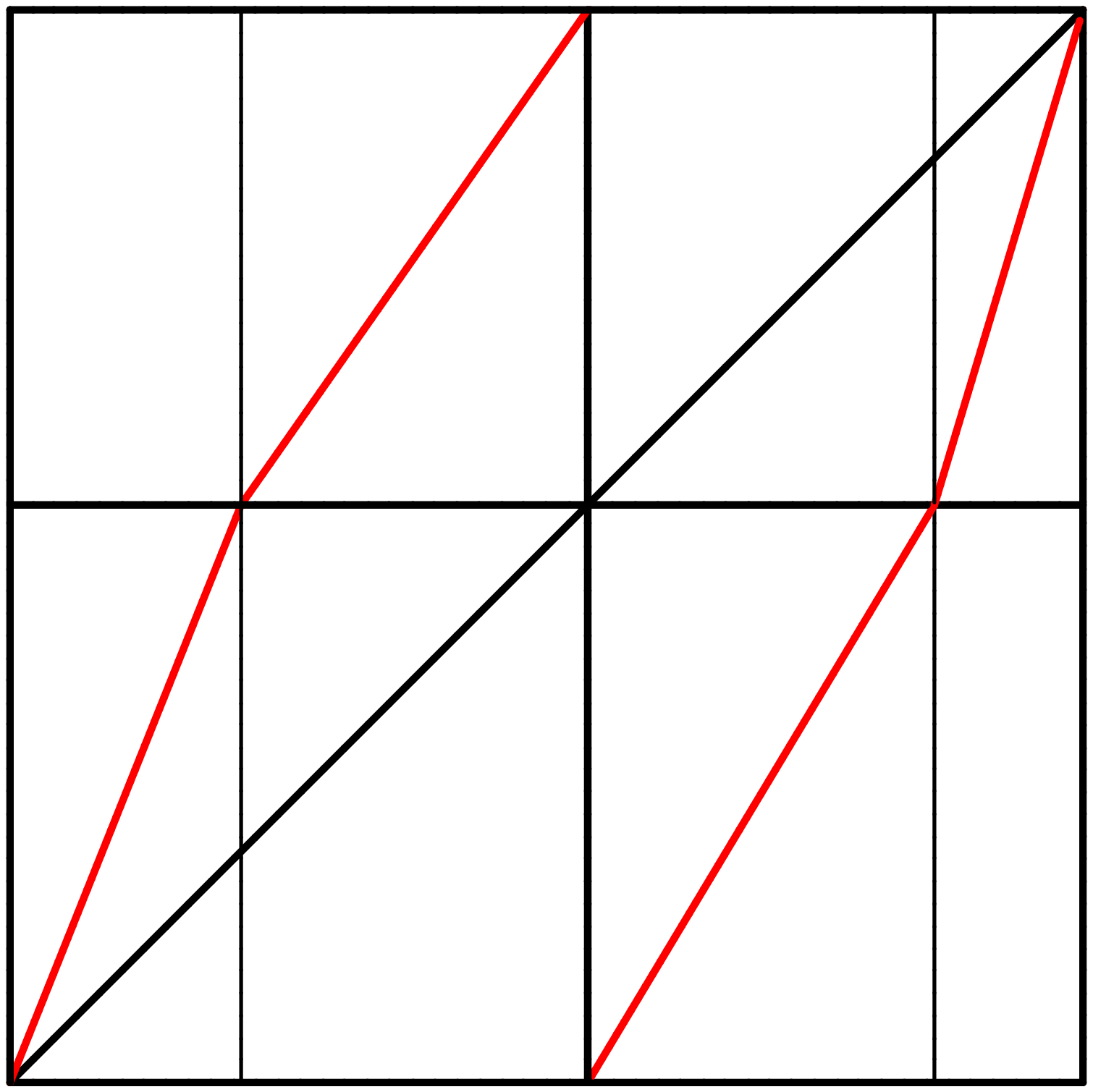}}

\put(30.9,15){\line(1,0){154}}
\put(30.9,10){\line(0,1){10}}
\put(84.8,10){\line(0,1){10}}
\put(184.9,10){\line(0,1){10}}

\put(55,2){$I_0$}
\put(135,2){$I_1$}

\put(223,15){\line(1,0){154}}
\put(223,10){\line(0,1){10}}
\put(305.9,10){\line(0,1){10}}
\put(377,10){\line(0,1){10}}

\put(260,2){$I_0$}
\put(340,2){$I_1$}

\end{picture}
\caption{\label{figBernoulliMarkov}
mappings generating a Bernoulli source and a stationary Markov chain of order $1$.
In both cases, Lebesgue measure is the initial one.}
\end{figure}

\vskip 10pt
The following classical examples often appear in the literature:
let $p\in ]0,1[$, $I_0=[0,1-p]$ and $I_1=]1-p,1]$.
Let $T:I\to I$ be the only function which maps linearly and increasingly $I_0$ and $I_1$ onto $I$ (see Figure~\ref{figBernoulliMarkov} when $p=0.65$, left side).
Then, starting from Lebesgue measure, the corresponding  probabilistic dynamical source is Bernoulli:
the $Y_n$ are i.i.d.~and $\Proba (Y_0=1)=p$.
In the same vein, if $T$ is the mapping drawn on the right side of
Figure~\ref{figBernoulliMarkov}, starting from Lebesgue measure, the $\{ 0,1\}$-valued process
$(Y_n)_{n\in\g N}$ is Markov and stationary, with transition matrix
\[\left(\begin{array}{cc}0.4&0.6\\0.7&0.3 \end{array}\right).\]
The assertions on both examples are consequences of Thales theorem.
These two basic examples are particular cases of Theorem~\ref{equivalenceProcessus}.

\subsubsection{Ordered subdivisions and ordered partitions of the interval}
\label{sssection:subdivisions}

\begin{Def}
A family $(I_w)_{w\in\rond W}$ of subintervals of $I$ indexed by all finite words is said to be an
\emph{$\rond A$-adic subdivision} of  $I$ whenever

(i) for any $w\in\rond W$, $I_w$ is the disjoint union of $I_{w\alpha}$, $\alpha\in\rond A$;

(ii) for any $v,w\in\rond W$, if $v<w$ for the alphabetical order, then
$$
\forall x\in I_v,~\forall y\in I_{w},~x<y.
$$
\end{Def}

\begin{Rem}
For any integer $p\geq 2$, the usual $p$-adic subdivision of $I$ is a particular case of $\rond A$-adic
subdivision for which $|\rond A| = p$ and $| I_w| = p^{-|w|}$ for any finite word $w\in\rond W$. For a general $\rond A$-adic subdivision, the intervals associated with two $k$-length words need not have the same length. 
\end{Rem}

The inclusion relations between the subintervals $I_w$ of an $\rond A$-adic subdivision are thus coded by
the prefix order in the complete $|\rond A|$-ary planar tree.
In particular, for any $w\in\rond W$ and for any cutset $\rond K$ of the complete $|\rond A|$-ary tree,
$$
I_w=\bigcup _{v\in\rond K}I_{wv}
$$
(this union is a disjoint one; see Section~\ref{trees} for a definition of a cutset).

We will use the following convention for $\rond A$-adic subdivisions:
we require the intervals $I_v$ to be open on the left side and closed on the right side, except the ones of the
form $I_{0^n}$ that are compact. Obviously, if $\mu$ is any probability measure on $\rond R=\rond A^{\N}$, there exists a
unique $\rond A$-adic subdivision of $I$ such that $|I_w|=\mu (w\rond R)$ for any $w\in\rond W$.

\vskip 5pt
Given an $\rond A$-adic subdivision of $I$, we extend the notation $I_w$ to right-infinite words by
$$
\forall r\in\rond R,~I_r=\bigcap _{{w\in\rond W}\atop{w{\rm ~prefix~of~}r}}I_w.
$$

\begin{Def}
A family $(I_v)_{v\in V}$ of subintervals of $I$ indexed by a totally ordered set $V$ is said to define an
\emph{ordered topological partition} of $I$ when

(i) $I=\bigcup _{v\in V}\cl (I_v)$,

(ii) for any $v,v'\in V$, $v\neq v'\Longrightarrow\intr (I_v)\cap\intr (I_{v'})=\emptyset$,

(iii) for any $v,v'\in V$,
$$
v\leq v'\Longrightarrow \forall x\in I_v,~\forall x'\in I_{v'},~x\leq x'
$$

where $\cl (I_v)$ and $\intr (I_v)$ stand respectively for the closure and the interior of $I_v$.
We will denote
$$
I=\ordPart {v\in V}I_v.
$$
\end{Def}

We will use the following fact:
if $I=\ordPart {v\in V}I_v=\ordPart {v\in V}J_v$ are two ordered topological partitions of $I$ indexed by the same denumerable ordered set $V$, then $I_v=J_v$ for any $v\in V$ as soon as $|I_v|=|J_v|$ for any $v\in V$.


\subsubsection{Definition of the mapping $T$}
\label{sectionDefinitionT}

Let $(U_n)_{n\geq 0}$ be a SVLMC, defined by its probabilized context tree $(\rond T,(q_c)_{c\in\rond C})$ and a
stationary\footnote{
Note that this construction can be made replacing $\pi$ by any probability measure on $\rond L$.} probability
measure $\pi$ on $\rond L$.
We first associate with $\pi$ the unique $\rond A$-adic subdivision $(I_w)_{w\in\rond W}$ of $I$, defined by:
$$
\forall w\in\rond W,~|I_w|=\pi (\overline w),
$$
(recall that if $w=\alpha _1\dots\alpha _N$, $\overline w$ is the reversed word $\alpha _N\dots\alpha _1$
and that $\pi (\overline w)$ denotes $\pi (\rond L\overline w)$).

\vskip 10pt
We consider now three ordered topological partitions of $I$.

\vskip 5pt
$\bullet$ The \emph{coding partition}

It consists in the family $(I_\alpha )_{\alpha\in\rond A}$:
$$
I=\ordPart{\alpha\in\rond A}I_\alpha =I_0+I_1.
$$

\vskip 5pt
$\bullet$ The \emph{vertical partition}

The denumerable set of finite and infinite contexts $\rond C$ is a cutset of the $\rond A$-ary tree.
The family $(I_c)_{c\in\rond C}$ thus defines the so-called vertical ordered topological partition
$$
I=\ordPart{c\in\rond C}I_c.
$$

\vskip 5pt
$\bullet$ The \emph{horizontal partition}

$\rond A\rond C$ is the set of leaves of the context tree
$\rond A\rond T=\{\alpha w,~\alpha\in\rond A,~w\in\rond T\}$.
As before, the family $(I_{\alpha c})_{\alpha c\in\rond A\rond C}$ defines the so-called horizontal ordered topological partition
$$
I=\ordPart{\alpha c\in\rond A\rond C}I_{\alpha c}.
$$
\begin{Def}
\label{map-on-context}
The mapping $T:I\to I$ is the unique left continuous function such that:

(i) the restriction of $T$ to any $I_{\alpha c}$ is affine and increasing;

(ii) for any $\alpha c\in\rond A\rond C$, $T(I_{\alpha c})=I_c$.
\end{Def}

The function $T$ is always increasing on $I_0$ and on $I_1$.
When $q_c(\alpha )\neq 0$, the slope of $T$ on an interval $I_{\alpha c}$ is
$1/q_c(\alpha )$. Indeed, with Formula~(\ref{probaStationnaireCalcul}), one has
$$|I_{\alpha c}|= \pi(\overline{c}\alpha) = q_c(\alpha)\pi(\overline{c}) = |I_c|q_c(\alpha ).$$
When $q_c(\alpha )=0$ and $|I_c|\neq 0$, the interval $I_{\alpha c}$ is empty so that $T$ is discontinuous at
$x_c=\pi (\{ s\in\rond L,~\overline s\leq c\})$ ($\leq$ denotes here the alphabetical order on $\rond R$).
Note that $|I_c|=0$ implies $|I_{\alpha c}|=0$.
In particular, if one assumes that all the probability measures $q_c$, $c\in\rond C$, are nontrivial
(\ie as soon as they satisfy $q_c(0)q_c(1)\neq 0$), then $T$ is continuous on $I_0$ and $I_1$.
Furthermore, $T(I_0)=\cl (T(I_1))=I$ and for any $c\in\rond C$, $T^{-1}I_c=I_{0c}\cup I_{1c}$
(see Figure~{\ref{figImageInverseIc}).

\begin{figure}[!h]
\begin{picture}(300,160)
\thicklines
\put(75,25){\includegraphics[width=270pt]{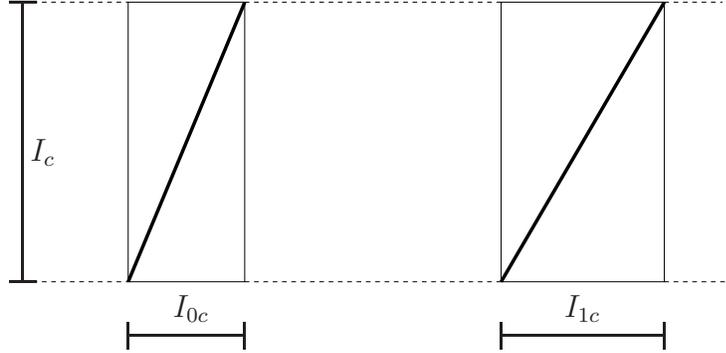}}
\put(132,13){$I_{0c}$}
\put(280,13){$I_{1c}$}
\put(115,6){\line(1,0){44}}
\put(115,1){\line(0,1){10}}
\put(159,1){\line(0,1){10}}
\put(256,6){\line(1,0){61.9}}
\put(256,1){\line(0,1){10}}
\put(317.9,1){\line(0,1){10}}

\put(75,26.2){\line(0,1){105.7}}
\put(78,72){$I_c$}
\put(70,26.3){\line(1,0){10}}
\put(70,132){\line(1,0){10}}
\end{picture}
\caption{\label{figImageInverseIc}
action of $T$ on horizontal and vertical partitions.
On this figure, $c$ is any context and the alphabet is $\{ 0,1\}$.}
\end{figure}

\vskip 5pt
{\bf Example: the four flower bamboo}

The \emph{four flower bamboo} is the VLMC defined by the finite probabilized context tree of
Figure~\ref{figArbreMapBambou4Fleurs}. There exists a unique stationary measure $\pi$  under conditions which are detailed later, in Example~\ref{bb4f}. We represent on Figure~\ref{figArbreMapBambou4Fleurs} the mapping $T$ built with this $\pi$, together with the respective subdivisions of $x$-axis and $y$-axis by the four $I_c$ and the eight $I_{\alpha c}$.
The $x$-axis is divided by both coding and horizontal partitions;
the $y$-axis is divided by both coding and vertical partitions.
This figure has been drawn with the following data on the four flower bamboo:
$q_{00}(0)=0.4$, $q_{010}(0)=0.6$, $q_{011}(0)=0.8$ and $q_1(0)=0.3$.

\begin{figure}[!h]
\begin{picture}(400,220)
\put(15,45){\includegraphics[height=120pt]{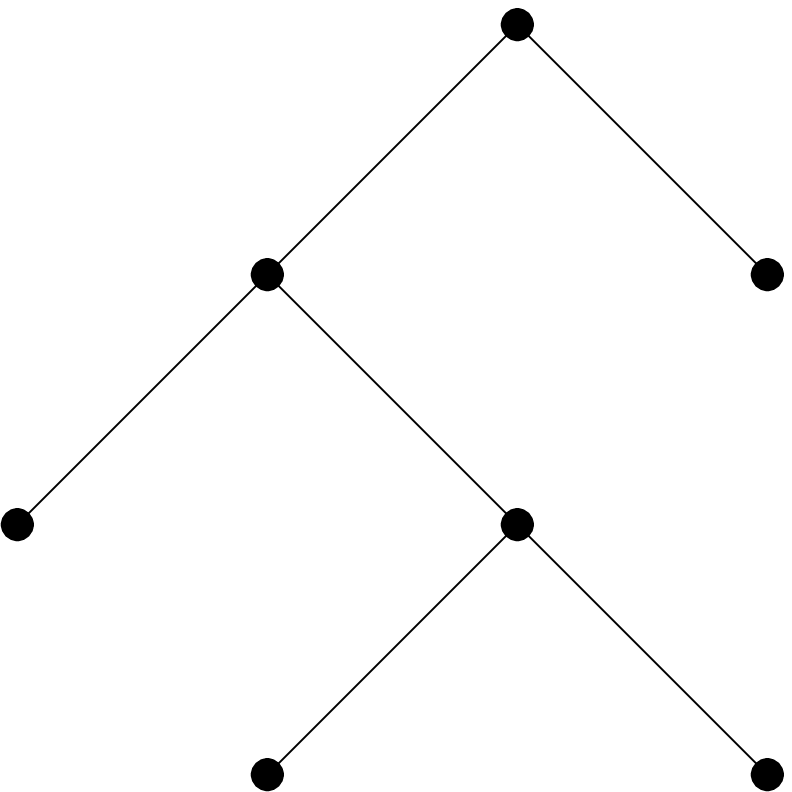}}
\put(225,30){\includegraphics[height=200pt]{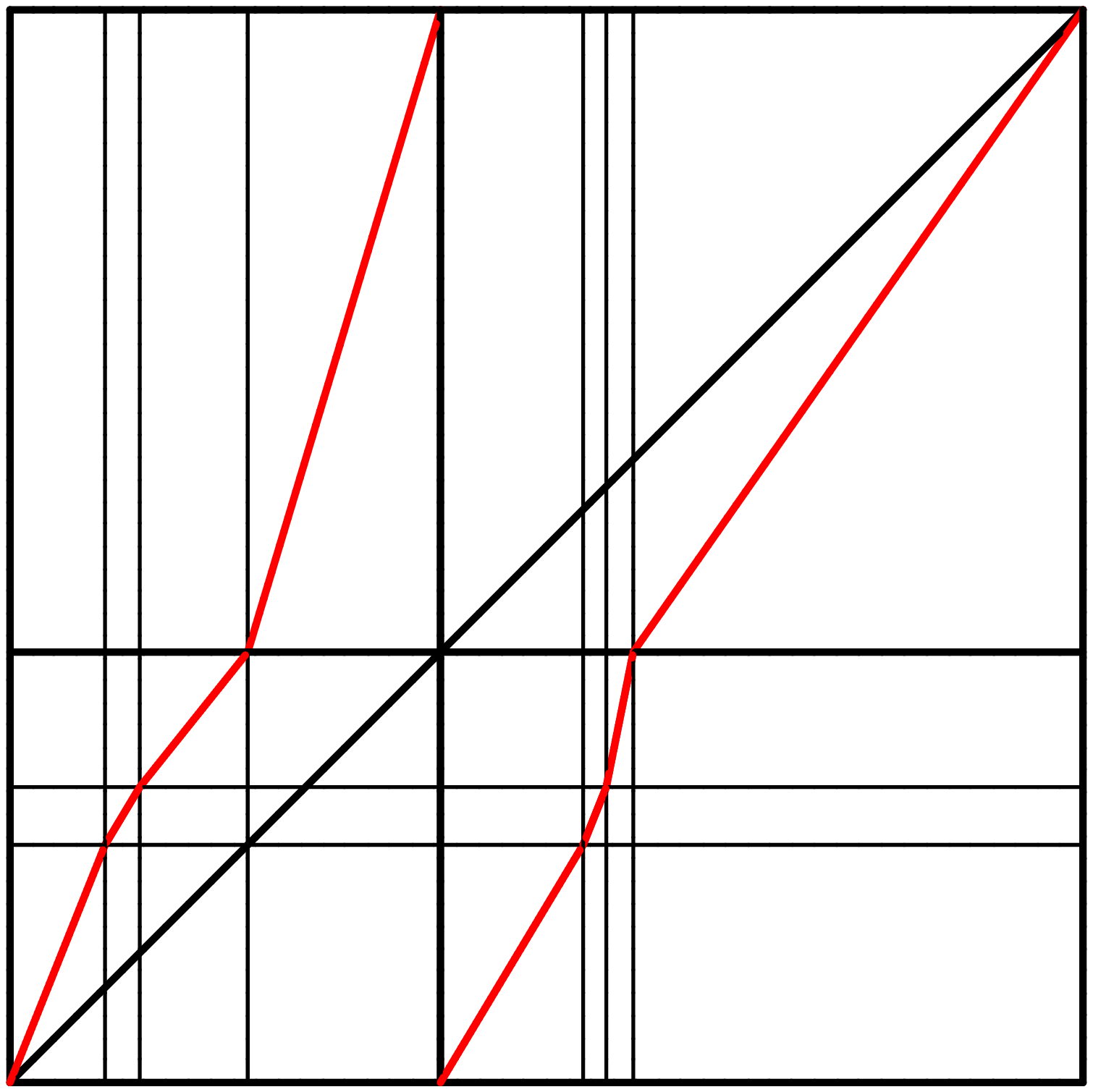}}
\put(45,30){$q_{010}$}
\put(124,30){$q_{011}$}
\put(11,67){$q_{00}$}
\put(129,105){$q_{1}$}

\thicklines
\put(248,37){\line(1,0){153.8}}
\put(275,41){$I_0$}
\put(355,41){$I_1$}
\put(248,34){\line(0,1){6}}
\put(309.8,34){\line(0,1){6}}
\put(401.8,34){\line(0,1){6}}

\put(248,17){\line(1,0){153.8}}
\put(245,23){$I_{000}$}
\put(289,23){$I_{01}$}
\put(310,23){$I_{100}$}
\put(360,23){$I_{11}$}
\put(245,-2){$I_{0010}$}
\put(275,-2){$I_{0011}$}
\put(310,-2){$I_{1010}$}
\put(342,-2){$I_{1011}$}
{\thinlines
\put(259,5){\vector(1,2){5}}
\put(289,5){\vector(-4,3){13}}
\put(324,3.8){\vector(2,3){7}}
\put(354,5){\vector(-2,1){17.5}}
}
\put(248,14){\line(0,1){6}}
{\thinlines
\put(261.7,14){\line(0,1){6}}
\put(266.7,14){\line(0,1){6}}
\put(282.1,14){\line(0,1){6}}
}
\put(309.8,14){\line(0,1){6}}
{\thinlines
\put(330.1,14){\line(0,1){6}}
\put(333.5,14){\line(0,1){6}}
\put(337.5,14){\line(0,1){6}}
}
\put(401.8,14){\line(0,1){6}}

\put(227,53){\line(0,1){153.9}}
\put(232,80){$I_0$}
\put(232,152){$I_1$}
\put(224,53){\line(1,0){6}}
\put(224,114.9){\line(1,0){6}}
\put(224,206.9){\line(1,0){6}}

\put(197,53){\line(0,1){153.9}}
\put(202,68){$I_{00}$}
\put(202,87){$I_{010}$}
\put(202,102){$I_{011}$}
\put(202,152){$I_1$}
\put(194,53){\line(1,0){6}}
\put(194,114.9){\line(1,0){6}}
\put(194,206.9){\line(1,0){6}}
{\thinlines
\put(194,87.2){\line(1,0){6}}
\put(194,95.4){\line(1,0){6}}
}
\end{picture}
\caption{\label{figArbreMapBambou4Fleurs}on the left, the 4 flower bamboo context tree.
On the right, its mapping together with the coding, the vertical and the horizontal partitions of $[0,1]$.}
\end{figure}

\subsubsection{Properties of the mapping $T$}
\label{sssec:proprietesT}

The following key lemma explains the action of the mapping $T$ on the intervals of the $\rond A$-adic subdivision $(I_w)_{w\in\rond W}$. More precisely, it extends the relation $T(I_{\alpha c})=I_c$, for any $\alpha c\in\rond A\rond C$, to any finite word. 
\begin{Lem}
\label{TIw}
For any finite word $w\in\rond W$ and any letter $\alpha\in\rond A$, $T(I_{\alpha w})=I_w$.
\end{Lem}

\pff
Assume first that $w\notin\rond T$. Let then $c\in\rond C$ be the unique context such that $c$ is a prefix of $w$.
Because of the prefix order structure of the $\rond A$-adic subdivision $(I_v)_v$, one has the first
ordered topological partition
\begin{equation}
\label{premierePartition}
I_c=\ordPart{{v\in\rond W,~|v|=|w|}\atop{c{\rm ~prefix~of~}v}}I_v
\end{equation}
(the set of indices is a cutset in the tree of  $c$ descendants).
On the other hand, the same topological partition applied to the finite word $\alpha w$ leads to
$$I_{\alpha c}=\ordPart{{v\in\rond W,~|v|=|w|}\atop{c{\rm ~prefix~of~}v}}I_{\alpha v}.$$
Taking the image by $T$, one gets the second ordered topological partition
\begin{equation}
\label{secondePartition}
I_c=\ordPart{{v\in\rond W,~|v|=|w|}\atop{c{\rm ~prefix~of~}v}}T(I_{\alpha v}).
\end{equation}
Now, if $c$ is a prefix of a finite word $v$, $I_{\alpha v}\subseteq I_{\alpha c}$ and the restriction of $T$ to
$I_{\alpha c}$ is affine. By Thales theorem, it comes
\[|T(I_{\alpha v})|=|I_{\alpha v}|.\frac{|I_c|}{|I_{\alpha c}|}.\]
Since $\pi$ is a stationary measure for the VLMC, 
\[|I_{\alpha c}|= \pi(\overline{c}\alpha) = q_c(\alpha)\pi(\overline{c}) = |I_c|q_c(\alpha ).\]
Furthermore, one has $\pi (\overline v\alpha )=q_c(\alpha )\pi (\overline v)$.
Finally, $|T(I_{\alpha v})|=|I_v|$.
Relations~(\ref{premierePartition}) and~(\ref{secondePartition}) are two ordered denumerable topological
partitions, the components with the same indices being of the same length: the partitions are necessarily the same.
In particular, because $w$ belongs to the set of indices, this implies that $T(I_{\alpha w})=I_w$.

Assume now that $w\in\rond T$.
Since the set of contexts having $w$ as a prefix is a cutset of the tree of the descendants of $w$, one has the disjoint union
$$
I_{\alpha w}=\bigcup _{{c\in\rond C}\atop{w{\rm ~prefix~of~}c}}I_{\alpha c}.
$$
Taking the image by $T$ leads to
$$T(I_{\alpha w})=\bigcup _{{c\in\rond C}, {w{\rm ~prefix~of~}c}}I_{c}=I_w$$
 and the proof is complete.
\QED

\begin{Rem}
The same proof shows in fact that if $w$ is any finite word, $T^{-1}I_w=I_{0w}\cup I_{1w}$ (disjoint union).
\end{Rem}

\begin{Lem}
\label{imageInverseBorelien}
For any $\alpha\in\rond A$, for any context $c\in\rond C$, for any Borelian set $B\subseteq I_c$,
$$
|I_\alpha\cap T^{-1}B|=|B|q_c(\alpha).
$$
\end{Lem}

\pff
It is sufficient to show the lemma when $B$ is an interval. The restriction of $T$ to $I_{\alpha c}$ is affine and $T^{-1}I_c=I_{0c}\cup I_{1c}$.
The result is thus due to Thales Theorem.
\QED

\begin{Cor}
\label{invarianceT}
If $T$ is the mapping associated with a SVLMC, Lebesgue measure is invariant by $T$,
\ie $|T¬^{-1}B|=|B|$ for any Borelian subset of $I$.
\end{Cor}

\pff
Since $B=\bigcup _{c\in\rond C}B\cap I_c$ (disjoint union), it suffices to prove that $|T¬^{-1}B|=|B|$ for
any Borelian subset of $I_c$ where $c$ is any context.
If $B\subseteq I_c$, because of Lemma~\ref{imageInverseBorelien},
\[|T^{-1}B|=|I_0\cap T^{-1}B|+|I_1\cap T^{-1}B|=|B|(q_c(0)+q_c(1))=|B|.\]
\QED

\subsubsection{SVLMC as dynamical source}
\label{sssection:VLMCisDS}

We now consider the stationary probabilistic dynamical source
$\left((I_\alpha )_{\alpha\in\rond A},\rho ,T, |.|\right)$ built from the SVLMC.
It provides the
$\rond A$-valued random process $(Y_n)_{n\in\N}$ defined by
$$
Y_n=\rho(T^n\xi )
$$
where $\xi$ is a uniformly distributed $I$-valued random variable and $\rho$ the coding function.
Since Lebesgue measure is $T$-invariant, all random variables $Y_n$ have the same law, namely
$\Proba (Y_n=0)=|I_0|=\pi (0)$.

\begin{Def}
Two $\rond A$-valued random processes $(V_n)_{n\in\g N}$ and $(W _n)_{n\in\g N}$ are called
\emph{symmetrically distributed} whenever for any $N\in\g N$ and for any finite word $w\in\rond A^{N+1}$,
$\Proba (W_0W_1\dots W_N=w)=\Proba (V_0V_1\dots V_N=\overline w)$.
\end{Def}

In other words, $(V_n)_{n\in\g N}$ and $(W_n)_{n\in\g N}$ are symmetrically distributed if and only if
for any $N\in\g N$, the random words $W_0W_1\dots W_N$ and
$V_NV_{N-1}\dots V_0$ have the
same distribution.

\begin{Th}
\label{equivalenceProcessus}
Let $(U_n)_{n\in\g N}$ be a SVLMC and $\pi$ a stationary probability measure on~$\rond L$.
Let $(X_n)_{n\in\g N}$ be the process of final letters of $(U_n)_{n\in\g N}$.
Let $T:I\to I$ be the mapping defined in Section~\ref{sectionDefinitionT}.
Then,

(i)
Lebesgue measure is $T$-invariant.

(ii)
If $\xi$ is any uniformly distributed random variable on $I$, the processes $(X_n)_{n\in\g N}$ and
$(\rho(T^n\xi ))_{n\in\g N}$ are symmetrically distributed.
\end{Th}

\vskip 10pt
\pff
\emph{(i)} has been already stated and proven in Corollary~\ref{invarianceT}.

\emph{(ii)} As before, for any finite word $w=\alpha _0\dots\alpha _N\in\rond W$, let
$B_w=\bigcap _{k=0}^NT^{-k}I_{\alpha _k}$ be the Borelian set of real numbers $x$ such that the right-infinite
sequence $(\rho (T^nx))_{n\in\g N}$ has $w$ as a prefix.
By definition, $B_\alpha =I_\alpha$ if $\alpha\in\rond A$.
More generally, we prove the following claim:
\emph{for any $w\in\rond W$, $B_w=I_w$.}
Indeed, if $\alpha\in\rond A$ and $w\in\rond W$, $B_{\alpha w}=I_\alpha\cap T^{-1}B_w$;
thus, by induction on the length of $w$, $B_{\alpha w}=I_\alpha\cap T^{-1}I_w=I_{\alpha w}$, the last
equality being due to Lemma \ref{TIw}.
There is now no difficulty in finishing the proof:
if $w\in\rond W$ is any finite word of length $N+1$, then
$\Proba (X_0\dots X_N=\overline w)=\pi (\overline w)=|I_w|$.
Thus, because of the claim,
$\Proba (X_0\dots X_N=\overline w)=|B_w|=\Proba (Y_0\dots Y_N=w)$.
This proves the result.
\QED

\section{Examples}
\label{sec:ex}

\subsection{The infinite comb}
\label{ssec:comb}

\subsubsection{Stationary probability measures}
\label{sssec:combSVLMC}

Consider the probabilized context tree given on the left side of Figure~\ref{figPeigneInfini}.
\begin{figure}[!h]
\begin{picture}(400,220)
\put(10,35){\includegraphics[height=150pt]{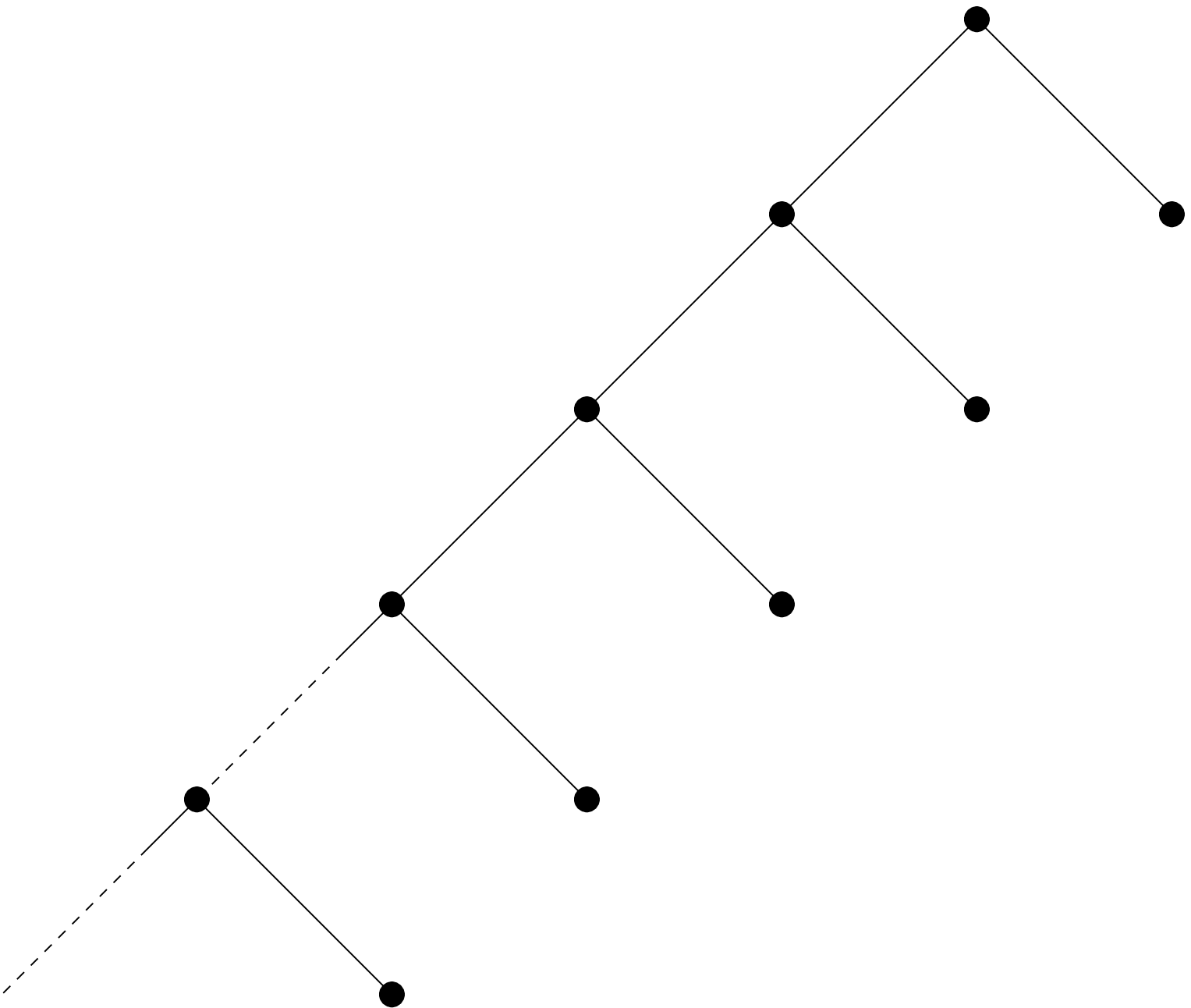}}
\put(220,10){\includegraphics[height=210pt]{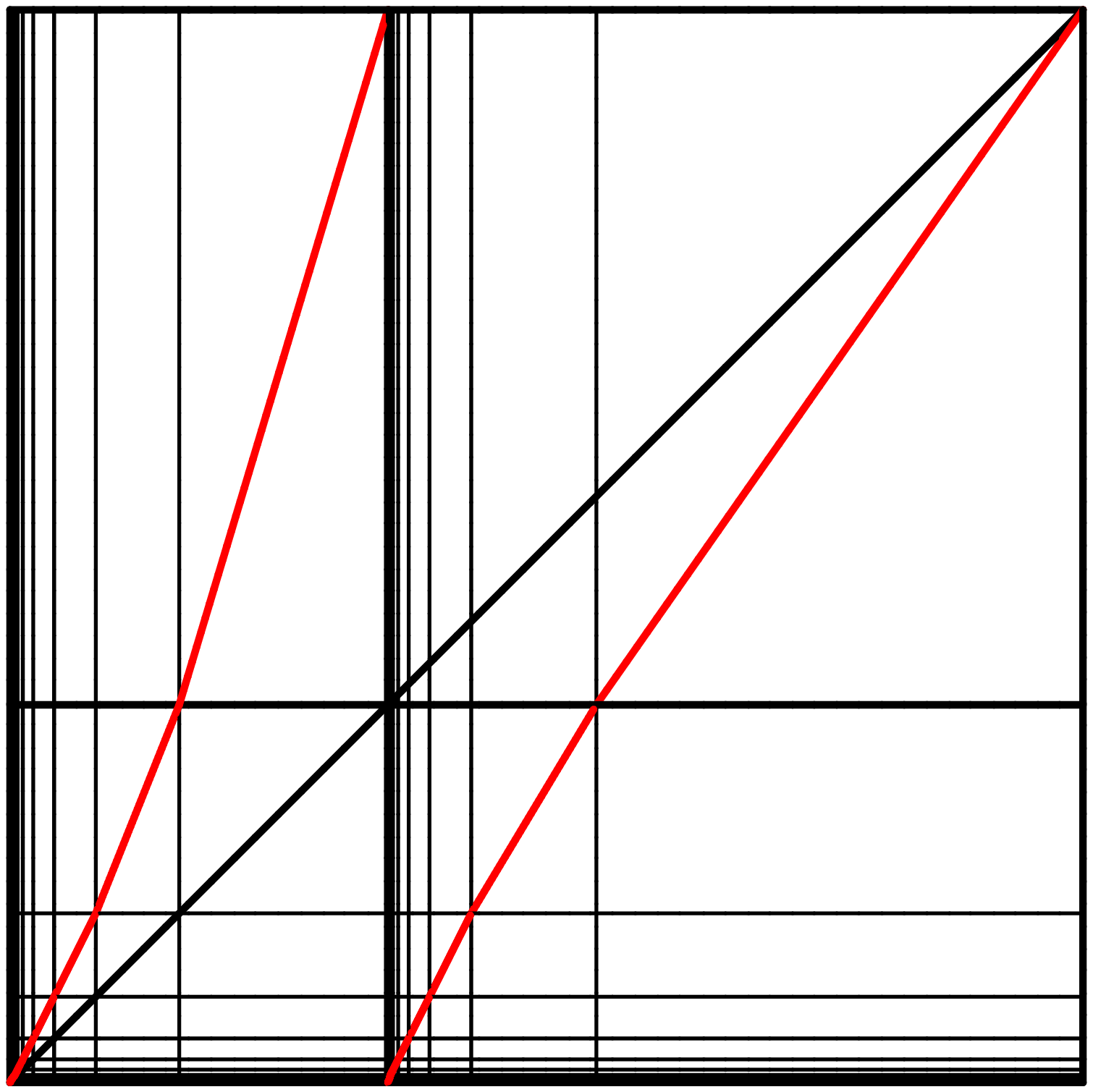}}

\thicklines
\put(4,20){$q_{0^\infty}$}
\put(60,20){$q_{0^n1}$}
\put(88,50){$q_{0001}$}
\put(117,80){$q_{001}$}
\put(150,110){$q_{01}$}
\put(180,140){$q_{1}$}

\put(244.4,13){\line(1,0){161.3}}
\put(270,20){$I_0$}
\put(346,20){$I_1$}
\put(244.4,10){\line(0,1){6}}
\put(301.5,10){\line(0,1){6}}
\put(405.7,10){\line(0,1){6}}

\put(220,34.2){\line(0,1){161.5}}
\put(225,60){$I_0$}
\put(225,135){$I_1$}
\put(217,34.2){\line(1,0){6}}
\put(217,90.9){\line(1,0){6}}
\put(217,195.7){\line(1,0){6}}
\end{picture}
\caption{\label{figPeigneInfini}
infinite comb probabilized context tree (on the left) and the associated dynamical system
(on the right).}
\end{figure}
In this case, there is one infinite leaf $0^\infty$ and countably many finite leaves $0^n1$, $n\in\N$.
The data of a corresponding VLMC consists thus in probability measures on $\rond A=\{ 0,1\}$:
$$
q_{0^\infty}
{\rm ~and~~}
q_{0^n1},~n\in\g N.
$$

\vskip 5pt
Suppose that $\pi$ is a stationary measure on $\rond L$.
We first compute $\pi (w)$ (notation (\ref{notationpiw})) as a function of $\pi (1)$ when $\overline w$ is any context or any internal node.
Because of Formula~(\ref{probaStationnaireCalcul}), $\pi (10)=\pi (1)q_1(0)$ and an immediate induction
shows that, for any $n\geq 0$,
\begin{equation}
\label{probaContextesPeigne}
\pi (10^n)=\pi (1)c_n,
\end{equation}
where $c_0=1$ and, for any $n\geq 1$,
\begin{equation}
\label{definitioncnPeigne}
c_n=\prod _{k=0}^{n-1}q_{0^k1}(0).
\end{equation}
The stationary probability of a reversed context is thus necessarily given by
Formula~(\ref{probaContextesPeigne}).
Now, if $0^n$ is any internal node of the context tree, we need going down along the branch in $\rond T$ to reach the contexts; using then the disjoint union $\pi (0^{n+1})=\pi (0^n)-\pi (10^n)$, by induction, it comes for any $n\geq 0$, 
\begin{equation}
\label{probaNoeudsInternesPeigne}
\pi (0^n)=1-\pi (1)\sum _{k=0}^{n-1}c_k.
\end{equation}
The stationary probability of a reversed internal node of the context tree is thus necessarily given by
Formula~(\ref{probaNoeudsInternesPeigne}).

\vskip 5pt
It remains to compute $\pi (1)$.
The denumerable partition of the whole probability space given by all cylinders based on leaves in the context
tree (Formula (\ref{partition})) implies
$1-\pi (0^\infty )=\pi (1)+\pi (10)+\pi (100)+\dots$, \ie
\begin{equation}
\label{partitionPeigne}
1-\pi (0^\infty )=\sum_{n\geq 0}\pi (1)c_n.
\end{equation}

This leads to the following statement that covers all cases of existence, unicity and nontriviality for a stationary
probability measure for the infinite comb.
In the generic case (named \emph{irreducible} case hereunder), we give a necessary and sufficient condition on the data for the existence of a stationary probability measure;
moreover, when a stationary probability exists, it is unique.
The \emph{reducible} case is much more singular and gives rise to nonunicity. 

\begin{Prop}
\label{mesureStationnairePeigne}
{\bf (Stationary probability measures for an infinite comb)}

Let $(U_n)_{n\geq 0}$ be a VLMC defined by a probabilized infinite comb.

\medskip

{(i) \emph{Irreducible case}}

\smallskip
Assume that $q_{0^\infty}(0)\neq 1$.

\medskip 
{(i.a) \emph{Existence}}

The Markov process $(U_n)_{n\geq 0}$ admits a stationary probability measure on $\rond L$ if and only if
the numerical series $\sum c_n$ defined by~(\ref{definitioncnPeigne}) converges.

\medskip
{(i.b) \emph{Unicity}}

Assume that the series $\sum c_n$ converges and denote
\begin{equation}
\label{S(1)Peigne}
S(1)=\sum _{n\geq 0}c_n.
\end{equation}

Then, the stationary probability measure $\pi$ on $\rond L$ is unique;
it is characterized by 
\begin{equation}
\label{pi(1)Peigne}
\pi (1)=\frac 1{S(1)}
\end{equation}
and Formulae~(\ref{probaContextesPeigne}),
(\ref{probaNoeudsInternesPeigne}) and~(\ref{probaStationnaireCylindres}).

Furthermore, $\pi$ is trivial if and only if $q_1(0)=0$, in which case it is defined by $\pi (1^\infty )=1$.

\bigskip

{(ii) \emph{Reducible case}}

\smallskip

Assume that $q_{0^\infty}(0)=1$.

\medskip 

{(ii.a)}
If the series $\sum c_n$ diverges, then the trivial probability measure $\pi$ on $\rond L$ defined by
$\pi (0^\infty )=1$ is the unique stationary probability.

\medskip

{(ii.b)}
If the series $\sum c_n$ converges, then there is a one parameter family of stationary probability measures
on $\rond L$.
More precisely, for any $a\in [0,1]$, there exists a unique stationary probability measure $\pi _a$ on $\rond L$
such that $\pi _a(0^\infty)=a$.
The probability $\pi _a$ is characterized by
$\pi _a(1)=\frac {1-a}{S(1)}$ and Formulae~(\ref{probaContextesPeigne}),
(\ref{probaNoeudsInternesPeigne}) and~(\ref{probaStationnaireCylindres}).

Furthermore, $\pi _a$ is non trivial except in the two following cases:\\
$\bullet$ $a=1$, in which case $\pi _1$ is defined by $\pi _1(0^\infty )=1$;\\
$\bullet$ $a=0$ and $q_1(0)=0$, in which case $\pi _0$ is defined by $\pi _0(1^\infty )=1$.
\end{Prop}

\pff
{\it (i)}
Assume that $q_{0^\infty}(0)\neq 1$ and that $\pi$ is a stationary probability measure.
By definition of probability transitions, $\pi (0^\infty )=\pi (0^\infty )q_{0^\infty}(0)$ so that $\pi (0^\infty )$
necessarily vanishes.
Thus, thanks to~(\ref{partitionPeigne}), $\pi (1)\neq 0$, the series $\sum c_n$ converges
and Formula~(\ref{pi(1)Peigne}) is valid.
Moreover, when $\overline w$ is any context or any internal node of the context tree, $\pi (w)$ is
necessarily given by Formulae~(\ref{pi(1)Peigne}), (\ref{probaContextesPeigne})
and~(\ref{probaNoeudsInternesPeigne}). This shows that for any finite word $w$, $\pi(w)$ is determined by Formula~(\ref{probaStationnaireCylindres}).
Since the cylinders $\rond L w$, $w\in\rond W$ span the $\sigma$-algebra on $\rond L$, there is at most one stationary
probability measure.
This proves the \emph{only if} part of {\it (i.a)}, the unicity and the characterization claimed in~{\it (i.b)}.

Reciprocally, when the series converges, Formulae~(\ref{pi(1)Peigne}), (\ref{probaContextesPeigne}),
(\ref{probaNoeudsInternesPeigne}) and~(\ref{probaStationnaireCylindres}) define a probability measure on
the semiring spanned by cylinders, which extends to a stationary probability measure on the whole
$\sigma$-algebra on $\rond L$
(see \cite{Billingsley} for a general treatment on semirings,
$\sigma$-algebra, definition and characterization of probability measures).
This proves the \emph{if} part of {\it (i.a)}.
Finally, the definition of $c_n$ directly implies that $S(1)=1$ if and only if $q_1(0)=0$.
This proves the assertion of {\it (i.b)} on the triviality of $\pi$.

\medskip

{\it (ii)}
Assume that $q_{0^\infty}(0)=1$.
Formula~(\ref{partitionPeigne}) is always valid so that the divergence of the series $\sum c_n$ forces $\pi (1)$ to vanish and, consequently, any stationary
measure $\pi$ to be the trivial one defined by $\pi (0^\infty )=1$.

Besides, with the assumption $q_{0^\infty}(0)=1$, one immediately sees that this trivial probability is stationary,
proving {\it (ii.a)}.

To prove {\it (ii.b)}, assume furthermore that the series $\sum c_n$ converges and let $a\in [0,1]$.
As before, any stationary probability measure $\pi$ is completely determined by $\pi (1)$.
Moreover, the probability measure defined by $\pi _a(1)=\frac {1-a}{S(1)}$,
Formulae~(\ref{probaContextesPeigne}),
(\ref{probaNoeudsInternesPeigne}) and~(\ref{probaStationnaireCylindres}) and standardly extended to the
whole $\sigma$-algebra on $\rond L$ is clearly stationary.
Because of Formula~(\ref{partitionPeigne}), it satisfies 
$$\pi _a(0^\infty )=1-\pi _a(1)S(1)=a.$$
This proves the assertion on the one parameter family.
Finally, $\pi _a$ is trivial only if $\pi _a(1)\in\{ 0,1\}$.
If $a=1$ then $\pi _a(1)=0$ thus $\pi _1$ is the trivial probability that only charges $0^\infty$.
If $a=0$ then $\pi _a(1)=1/S(1)$ is nonzero and it equals $1$ if and only if $S(1)=1$, \ie if and only if
$q_1(0)=0$, in which case $\pi _0$ is the trivial probability that only charges $1^\infty$.
\QED

\begin{Rem}
This proposition completes previous results which give sufficient conditions for the existence of a stationary measure for an infinite comb.
For instance, in \cite{Galves-Locherbach}, the intervening condition is 
$$ \sum_{k\geq 0}q_{0^k1}(1)=+\infty,$$
which is equivalent with our notations to $c_n \to 0$. Note that if $\sum c_n$ is divergent, then the only possible stationary distribution is the trivial Dirac measure $\delta_{0^\infty}$. 
\end{Rem}

\subsubsection{The associated dynamical system}
\label{sssec:combSD}

The vertical partition is made of the intervals $I_{0^{n}1}$ for $n\geq 0$. The horizontal partition consists in the 
intervals $I_{00^{n}1}$ and $I_{10^{n}1}$, for $n\geq 0$, together with two intervals coming from the infinite context, namely
$I_{0^{\infty}}$ and $I_{10^{\infty}}$. In the irrreducible case, $\pi(0^{\infty})=0$ and
these two last intervals become two accumulation points of the partition: $0$ and $\pi(0)$. The following lemma is classical and helps understand the behaviour of the mapping $T$ at these accumulation points.
\begin{Lem}\label{TAF}
Let $f:[a,b]\to{\R}$ be continuous on $[a,b]$, differentiable on $]a,b[ \setminus D$ where $D$ is a countable set. The fonction $f$ admits a right derivative at $a$
and
$$
f_r'(a)=\lim_{{x\to a,x>a}\atop{x\notin D}} f'(x)
$$
as soon as this limit exists.
\end{Lem}

\begin{Cor}
If $(q_{0^n1}(0))_{n\in{\N}}$ converges, then $T$ is differentiable at $0$ and
$\pi (0)$ (with a possibly infinite derivative) and 
 $$
 T'_r(0)=\lim_{n\to +\infty}\frac{1}{q_{0^n1}(0)},\ \ T'_r(\pi (0))=\lim_{n\to +\infty}\frac{1}{q_{0^n1}(1)}.
 $$
\end{Cor}
When $(q_{0^n1}(0))_{n\in{\N}}$ converges to $1$, $T'_r(0)=1$. In this case, $0$ is an indifferent fixed point and $T'_r(\pi (0)) = +\infty$. The mapping $T$ is a slight modification of the so-called Wang map (\cite{Wang}). The statistical properties of the Wang map are quite well understood (\cite{Lambert}). The corresponding dynamical source is said intermittent.

\subsubsection{Dirichlet series.}
For a stationary infinite comb, the Dirichlet series is defined on a suitable vertical open strip of $\g C$ as
$$\Lambda (s)=\sum _{w\in\rond W}\pi (w)^s.$$
In the whole section we suppose that $\sum c_n$ is convergent. Indeed, if it is divergent then the only stationary measure is the Dirac measure $\delta_{0^\infty}$ and $\Lambda (s)$ is never defined.

\medskip

The computation of the Dirichlet series is tractable because of the following formula:
for any finite words $w,w'\in\rond W$,
\begin{equation}
\label{w1wPeigne}
\pi (w1w')\pi(1)=\pi (w1)\pi (1w').
\end{equation}
This formula, which comes directly from Formula~(\ref{probaStationnaireCylindres}), is true because of the very particular form of the contexts in the infinite comb. It is the expression of its renewal property.
The computation of the Dirichlet series is made in two steps.

\medskip
{\bf Step 1.}
A finite word either does not contain any $1$ or is of the form $w10^n$, $w\in\rond W$, $n\geq 0$.
Thus,
$$\Lambda (s)=
\sum _{n\geq 0}\pi (0^n)^s
+\sum _{n\geq 0}\sum _{w\in\rond W}\pi (w10^n)^s.$$
Because of Formulae~(\ref{w1wPeigne}) and~(\ref{pi(1)Peigne}),
$\pi (w10^n)=S(1)\pi (w1)\pi (10^n)$.
Let us denote
$$
\Lambda _1(s)=\sum _{w\in\rond W}\pi (w1)^s.
$$
With this notation and Formulae~(\ref{probaContextesPeigne}) and~(\ref{probaNoeudsInternesPeigne}),
$$
\Lambda (s)=\frac 1{S(1)^s}\sum _{n\geq 0}R_n^s
+\Lambda _1(s)\sum _{n\geq 0}c_n^s
$$
where $R_n$ stands for the rest
\begin{equation}
\label{definitionRnPeigne}
R_n=\sum _{k\geq n}c_k.
\end{equation}

{\bf Step 2.}
It consists in the computation of $\Lambda _1$.
A finite word having $1$ as last letter either can be written $0^n1$, $n\geq 0$ or is of the form $w10^n1$,
$w\in\rond W$, $n\geq 0$.
Thus it comes,
$$\Lambda _1(s)=\sum _{n\geq 0}\pi (0^n1)^s+\sum _{n\geq 0}\sum _{w\in\rond W}\pi (w10^n1)^s.$$
By Formulae~(\ref{w1wPeigne}) and~(\ref{probaContextesPeigne}),
$\pi (w10^n1)
=\pi (w1)c_nq_{0^n1}(1)
=\pi (w1)(c_n-c_{n+1})$, so that
$$\Lambda _1(s)=
\frac 1{S(1)^s}\sum _{n\geq 0}c_n^s
+\Lambda _1(s)\sum _{n\geq 0}(c_n-c_{n+1})^s$$
and
$$
\Lambda _1(s)=
\frac{1}{S(1)^s}
\cdot\frac{\sum _{n\geq 0}c_n^s}{1-\sum _{n\geq 0}(c_n-c_{n+1})^s}.
$$
Putting results of both steps together, we obtain the following proposition.

\begin{Prop}
With notations~(\ref{definitioncnPeigne}), (\ref{S(1)Peigne}) and~(\ref{definitionRnPeigne}), the Dirichlet series of
a source obtained from a stationary infinite comb is
$$
\Lambda (s)=
\frac{1}{S(1)^s}
\left[
\sum _{n\geq 0}R_n^s
+\frac{\left(\sum _{n\geq 0}c_n^s\right)^2}{1-\sum _{n\geq 0}(c_n-c_{n+1})^s}
\right].
$$
\end{Prop}

\begin{Rem}
The analytic function
$\frac{\left(\sum _{n\geq 0}c_n^s\right)^2}{1-\sum _{n\geq 0}(c_n-c_{n+1})^s}$
is always singular for $s=1$ because its denominator vanishes while its
numerator is a convergent series.
\end{Rem}

\vskip 5pt
{\bf Examples.}
 (1) Suppose that $0<a<1$ and that $q_{0^n1}(0)=a$ for any $n\geq 0$.
Then $c_n=a ^n$, $R_n=\frac{a ^n}{1-a}$ and $S(1)=\frac 1{1-a}$.
For such a source, the Dirichlet series is
$$
\Lambda (s)=\frac 1{1-[a ^s+(1-a)^s]}.
$$
In this case, the source is memoryless: all letters are drawn independently with the same distribution.
The Dirichlet series of such sources have been extensively studied in \cite{FRV} in the realm
of asymptotics of average parameters of a trie.

\vskip 5pt
(2) Extension of Example 1:
take $a,b\in ]0,1[$ and consider the probabilized infinite comb defined by
$$
q_{0^n1}(0)=\left\{\begin{array}{l}a~{\rm if}~n~{\rm is~even},\\b~{\rm if}~n~{\rm is~odd.}\end{array}\right.
$$
After computation, the Dirichlet series of the corresponding source under the stationary distribution turns out to have the explicit form
$$
\Lambda (s)=\frac{1}{1-(ab)^s}
\left[
1+\left(\frac{a+ab}{1+a}\right) ^s+\left(\frac{1-ab}{1+a}\right) ^s\frac {(1+a^s)^2}{1-(ab)^s-(1-a)^s-a^s(1-b)^s}
\right] .
$$
The configuration of poles of $\Lambda$ depends on arithmetic properties (approximation by rationals) of the logarithms of $ab$, $1-a$ and $a(1-b)$.
The poles of such a series are the same as in the case of a memoryless source with an alphabet of
three letters, see \cite{FRV}. This could be extended to a family of examples.

\vskip 5pt
(3) Let $\alpha >2$.
We take data $q_{0^n1}(0)$, $n\geq 0$ in such a way that $c_0=1$ and, for any
$n\geq 1$,
$$
c_n=\zeta (n,\alpha ):=\frac 1{\zeta (\alpha )}\sum _{k\geq n}\frac 1{k^\alpha},
$$
where $\zeta$ is the Riemann function. Since $c_n\in \mathcal{O}(n^{1-\alpha})$ when $n$ tends to infinity, there exists a unique stationary probability measure
$\pi$ on $\rond L$.
One obtains 
$$S(1)=1+\frac{\zeta (\alpha -1)}{\zeta (\alpha )}$$ and, for any $n\geq 1$,
$$R_n=\frac{\zeta (\alpha -1)}{\zeta (\alpha )}\zeta (n,\alpha -1)-(n-1)\zeta (n,\alpha ).$$
In particular, $R_n\in \mathcal{O}(n^{2-\alpha})$ when $n$ tends to infinity.
The final formula for the Dirichlet series of this source is
$$
\Lambda (s)=\frac{1}{S(1)^s}\left[\sum _{n\geq 0}R_n^s+\frac{\left(\sum _{n\geq 0}c_n^s\right)^2}{1-\frac{\zeta (\alpha s)}{\zeta (\alpha )^s}}\right].
$$

\vskip 5pt
(4) One case of interest is when the associated dynamical system has an indifferent fixed point
(see Section \ref{sssec:combSD}), for example when
$$
q_{0^n1}(0)=\left(1-\frac{1}{n+2}\right)^\alpha,
$$
with $1<\alpha <2$.
In this situation, $c_n=(1+n)^{-\alpha}$ and 
$$
\Lambda (s)=\sum _{n\geq 1}\zeta (n,\alpha )^s
+\frac{\zeta (\alpha s)^2}{\zeta (\alpha )^s}\cdot
\frac{1}{\displaystyle 1-\sum _{n\geq 1}\frac 1{n^{\alpha s}}\left[ 1-\left( 1-\frac 1{n+1}\right) ^\alpha\right]}.
$$

\subsubsection{Generating function for the exact distribution of word occurrences in a sequence generated by a comb}
\label{ssec:generatrice-peigne}

The behaviour of the entrance time into cylinders is a natural question arising in dynamical systems. There exists a large literature on the asymptotic properties of entrance times into cylinders for various kind of systems, symbolic or geometric; see \cite{AG} for an extensive review on the subject. Most of the results deal with an exponential approximation of the distribution of the first entrance time into a small cylinder, sometimes with error terms. The most up-to-date result on this framework is \cite{AS}, not published yet, in which the hypothesis are made only in terms of the mixing type of the source (so-called $\alpha$-mixing). We are here interested in exact distribution results instead of asymptotic behaviours. 

Several studies in probabilities on words are based on generating functions. For example one may cite \cite{Regnier}, \cite{Schbath}, \cite{SP}. 
For i.i.d.~sequences, \cite{Blom} give the generating function of the first occurrence of a word, based on a recurrence relation on the probabilities. This result is extended to Markovian sequences by \cite{Daudin-Robin}. Nonetheless, other approaches are considered: one of the more general techniques is the so-called Markov chain embedding method introduced by \cite{Fu} and further developped by \cite{Fu-Koutras}, \cite{Koutras}. A martingale approach (see \cite{Gerber-Li}, \cite{Li}, \cite{Williams}) is an alternative to the Markov chain embedding method. These two approaches are compared in \cite{Pozdnyakov}. \\

We establish results on the exact distribution of word occurrences in a random sequence generated by a comb (or a bamboo in Section \ref{ssec:generatrice-bambou}). More precisely, we make explicit  the generating function of the random variable giving the $r$\textsuperscript{th} occurrence of a $k$-length word, for any word $w$ such that $\overline{w}$ is not an internal node of $\rond T$.  \\
Let us consider the process $X=(X_n)_{n\geq 0}$ of final letters of $(U_n)_{n\geq 0}$, in the particular case of a SVLMC defined by an infinite comb. Let $w=w_1 \ldots w_k$ be a word of length $k\geq 1$. We say that $w$ occurs at position $n\geq k$ in the sequence $X$ if the word $w$ ends at position $n$:
\[\{w \mbox{ at }n\}=\{X_{n-k+1}\ldots X_{n}=w\}=\{U_n \in \rond Lw\}.\]
Let us denote by $T_w^{(r)}$ the position of the $r$\textsuperscript{th} occurrence of $w$ in $X$ and $\Phi_w^{(r)}$ its generating function:
$$
\Phi_w^{(r)} (x) := \sum_{n\geq 0} \Proba\left( T_w^{(r)} =n\right) x^n.
$$
The following notation is used in the sequel: for any finite word $u \in \rond W$, for any	finite context $c \in\rond C$ and for any $n\geq 0$,
\[q_c^{(n)}(u)=\Proba\left(X_{n-|u|+1}\ldots X_{n}=u|X_{-(|c|-1)}\ldots X_{0}=\overline c\right).\]
These quantities may be computed in terms of the data $q_c$. Proposition \ref{propgengen} generalizes results of \cite{Daudin-Robin}.
\begin{Prop}
\label{propgengen} 
For a SVLMC defined by an infinite comb, with the above notations,  for a word $w$ such that $\overline{w}$ is non internal node,  the generating function of its first occurrence is given, for $|x|<1$, by
$$
\Phi_w^{(1)}(x)=\frac{x^k\pi(w)}{(1-x)S_w(x)}
$$
and the generating function of its
$r$\textsuperscript{th} occurrence is given, for $|x|<1$, by
\[\Phi_w^{(r)}(x)=\Phi_w^{(1)}(x)\left(1-\frac{1}{S_w(x)}\right)^{r-1},\]
where
\begin{eqnarray*}
S_w(x)&=& C_w(x) + \sum_{j=k}^{\infty}q_{\petitlpref(w)}^{(j)}(w)x^j,\\
C_w(x)&=& 1+\sum_{j=1}^{k-1}\ind{w_{j+1}\ldots w_{k}=w_{1}\ldots w_{k-j}}q_{\petitlpref(w)}^{(j)}\left(w_{k-j+1} \ldots w_k\right)x^j.
\end{eqnarray*}
\end{Prop}
\begin{rem} The  term $C_w(x)$ is a generalization of the probabilized autocorrelation polynomial defined in \cite{Jacquet-Szpankowski} in  the particular case when the $(X_n)_{n\geq 0}$ are independent and identically distributed.  For a word $w=w_1 \ldots w_k$ this polynomial is equal to
\[c_w(x)=\sum_{j=0}^{k-1}c_{j,w}\frac 1{\pi(w_{1}\ldots w_{k-j})} x^j,\]
where $c_{j,w}=1$ if the $k-j$-length suffix of $w$ is equal to its $k-j$-length prefix, and is equal to zero otherwise. When the $(X_n)_{n\geq 0}$ are independent and identically distributed, we have
\begin{equation*}
\sum_{j=1}^{k-1}\ind{w_{j+1}\ldots w_{k}=w_{1}\ldots w_{k-j}}q_{\petitlpref(w)}^{(j)}\left(w_{k-j+1} \ldots w_k\right)x^j=\sum_{j=1}^{k-1}c_{j,w}\frac{\pi(w)}{
\pi\left(w_{1} \ldots w_{k-j}\right)}x^j
\end{equation*}
that is $$C_w(x)= \pi(w)c_w(x).$$
\end{rem}
\pff
We first deal with $w=10^{k-1}$, that is the only word $w$ of length $k$ such that $\overline w \in \rond C$. For the sake of shortness, we will denote by $p_n$ the probability that $T_w^{(1)}=n$. From the obvious decomposition
\[\{w \mbox{ at }n\}=\{T_w^{(1)}=n\}\cup\{T_w^{(1)}<n \mbox{ and }w \mbox{ at }n\}, \quad \mbox{(disjoint union)}\]
it comes by stationarity of $\pi$
\[\pi(w)=p_n+\sum_{z=k}^{n-1}p_z\Proba\left(X_{n-k+1}\ldots X_{n}=w|T_w^{(1)}=z\right).\]
Due to the renewal property of the comb, the conditional probability can be rewritten
\[\left\{
\begin{array}{ll}
 \Proba\left(X_{n-k+1}\ldots X_n=w|X_{z-k+1}\ldots X_z=w\right) & \mbox{if }z\leq n-k\\
 0 &  \mbox{if }z> n-k\\
\end{array}
\right.,
\]
the second equality is due to the lack of possible auto-recovering in $w$. Consequently, we have 
\[\pi(w)=p_n+\sum_{z=k}^{n-k}p_zq_{\overline w}^{(n-z)}(w).\]
Hence, for $x<1$, it comes
\[\frac{x^k\pi(w)}{1-x}=\sum_{n=k}^{+\infty}p_nx^n+\sum_{n=k}^{+\infty}x^n\sum_{z=k}^{n-k}p_zq_{\overline w}^{(n-z)}(w),\]
so that
\[\frac{x^k\pi(w)}{1-x}=\Phi_w^{(1)}(x)\left(1+\sum_{j=k}^{\infty}x^jq_{\overline w}^{(j)}(w)\right),\]
which leads to
\[\Phi_w^{(1)}(x)=\frac{x^k\pi(w)}{(1-x)S_w(x)}.\]
Note that when $w=10^{k-1}$, $C_w(x) = 1$.

Proceeding in the same way for the $r$\textsuperscript{th} occurrence, from the decomposition
\[\{w \mbox{ at }n\}=\{T_w^{(1)}=n\}\cup\{T_w^{(2)}=n\}\cup\ldots\cup\{T_w^{(r)}=n\}\cup\{T_w^{(r)}<n \mbox{ and }w \mbox{ at }n\},\]
and denoting $p(n,\ell)=\Proba(T_w^{(\ell)}=n)$, the following recursive equation holds:
\[\pi(w)=p_n+p(n,2)+\ldots + p(n,r)+\sum_{z=k}^{n-1}\Proba\left(T_w^{(r)}=z\mbox{ and }w \mbox{ at }n\right).\]
Again, by splitting the last term into two terms and using the non-overlapping structure of $w$, one gets
\[\pi(w)=p_n+p(n,2)+\ldots +p(n,r)+\sum_{z=k}^{n-k}p_zq_{\overline w}^{(n-z)}(w).\]
From this recursive equation, proceeding exactly in the same way, one gets for the generating function, for $x<1$,
\[\Phi_w^{(r)}(x)=\Phi_w^{(1)}(x)\left(1- \frac{1}{S_w(x)}\right)^{r-1}.
\]
Let us now consider the case of words $w$ such that $\overline w\notin \rond T$, that is the words $w$ such that $w_j=1$ for at least one integer $j\in \{2, \ldots, k\}$. We denote by $i$ the last position of a $1$ in $w$, that is $\lpref(w)=0^{k-i}1$. Once again we have
\[
\pi(w)=p_n+\sum_{z=k}^{n-1}p_z\Proba\left(X_{n-k+1}\ldots X_{n}=w|T_w^{(1)}=z\right).\]
When $z\leq n-k$, due to the renewal property, the conditional probability can be rewritten as
\[\Proba\left(X_{n-k+1}\ldots X_{n}=w|T_w^{(1)}=z\right)=q_{\petitlpref(w)}^{(n-z)}(w).\]

\begin{picture}(200,110)
\linethickness{1pt}
\put(20,62){\line(1,0){330}}
\put(330,59){\line(0,1){6}}
\put(327,69){$n$}
\put(200,59){\line(0,1){6}}
\put(197,69){$z$}
\linethickness{2pt}
\put(170,82){\line(1,0){160}}
\put(170,79){\line(0,1){6}}
\put(165,90){$w_1$}
\put(200,79){\line(0,1){6}}
\put(193,90){$w_{k-n+z}$}
\put(330,79){\line(0,1){6}}
\put(323,90){$w_k$}
\put(200,42){\line(-1,0){160}}
\put(40,39){\line(0,1){6}}
\put(37,30){$w_1$}
\put(170,39){\line(0,1){6}}
\put(152,30){$w_{n-z+1}$}
\put(200,39){\line(0,1){6}}
\put(197,30){$w_k$}
\end{picture}

When $z>n-k$ (see figure above), 
\begin{equation*}
\Proba\left(w \mbox{ at }n|T_w^1=z\right)=\ind{w_{n-z+1}\ldots w_k=w_{1}\ldots w_{k-n+z}}q_{\petitlpref(w)}^{(n-z)}(w_{k-n+z+1}\ldots w_k),
\end{equation*}
\begin{picture}(400,50)
\linethickness{1pt}
\put(320,22){$w=*10\cdots 0$}
\put(357,32){$\overbrace{\rule{1cm}{0cm}}$}
\put(364,42){$\scriptstyle k-i$}
\put(353.5,7){\vector(0,1){10}}
\put(340,0){$\scriptstyle n-k+i$}
\put(383.5,7){\vector(0,1){10}}
\put(382,0){$\scriptstyle n$}
\put(0,15){\begin{minipage}{100 truemm}
this equality holding if $n-k+i\neq z$.
But when $z=n-k+i$, because the first occurrence of $w$ is at $z$, necessarily $w_k=1$ and hence $i=k$, and $z=n$ which contradicts $z<n$.
Consequently for $z=n-k+i$ we have
\end{minipage}}
\end{picture}
$$
\Proba\left(X_{n-k+1}\ldots X_{n}=w|T_w^1=z\right)=0 = \ind{w_{n-z+1}\ldots w_k=w_{1}\ldots w_{k-n+z}}.
$$
Finally one gets
\begin{eqnarray*}
\pi(w)=p_n &+& \sum_{z=1}^{n-k}p_zq_{\petitlpref(w)}^{(n-z)}(w)\\
&+&\sum_{z=n-k+1}^{n-1}p_z\ind{w_{n-z+1}\ldots w_k=w_1\ldots w_{k-n+z}}q_{\petitlpref(w)}^{(n-z)}(w_{k-n+z+1}\ldots w_k),
\end{eqnarray*}
and hence
\[\Phi_w^{(1)}(x)=\frac{x^k\pi(w)/(1-x)}{\displaystyle 1+\sum_{j=k}^{\infty}x^jq_{\petitlpref(w)}^{(j)}(w)+\sum_{j=1}^{k-1}x^j\ind{w_{j+1}\ldots w_{k}=w_1...w_{k-j}}q_{\petitlpref(w)}^{(j)}(w_{k-j+1}\ldots w_k)}.\]
Proceeding exactly in the same way by induction on $r$, we get the expression of Theorem~\ref{propgengen} for the $r$-th occurrence.
\QED

\begin{Rem}
The case of internal nodes $w=0^k$ is more intricate, due to the absence of any symbol $1$ allowing a renewal argument. Nevertheless, for the forthcoming applications, we will not need the explicit expression of the generating function of such words occurrences.
\end{Rem}

\subsection{The Bamboo blossom}
\label{ssec:blossom}
\subsubsection{Stationary probability measures}
\label{sssec:blossomSVLMC}

Consider the probabilized context tree given by the left side of Figure~\ref{figBambou}.
\begin{figure}[!h]
\begin{picture}(400,210)
\put(40,10){\includegraphics[height=190pt]{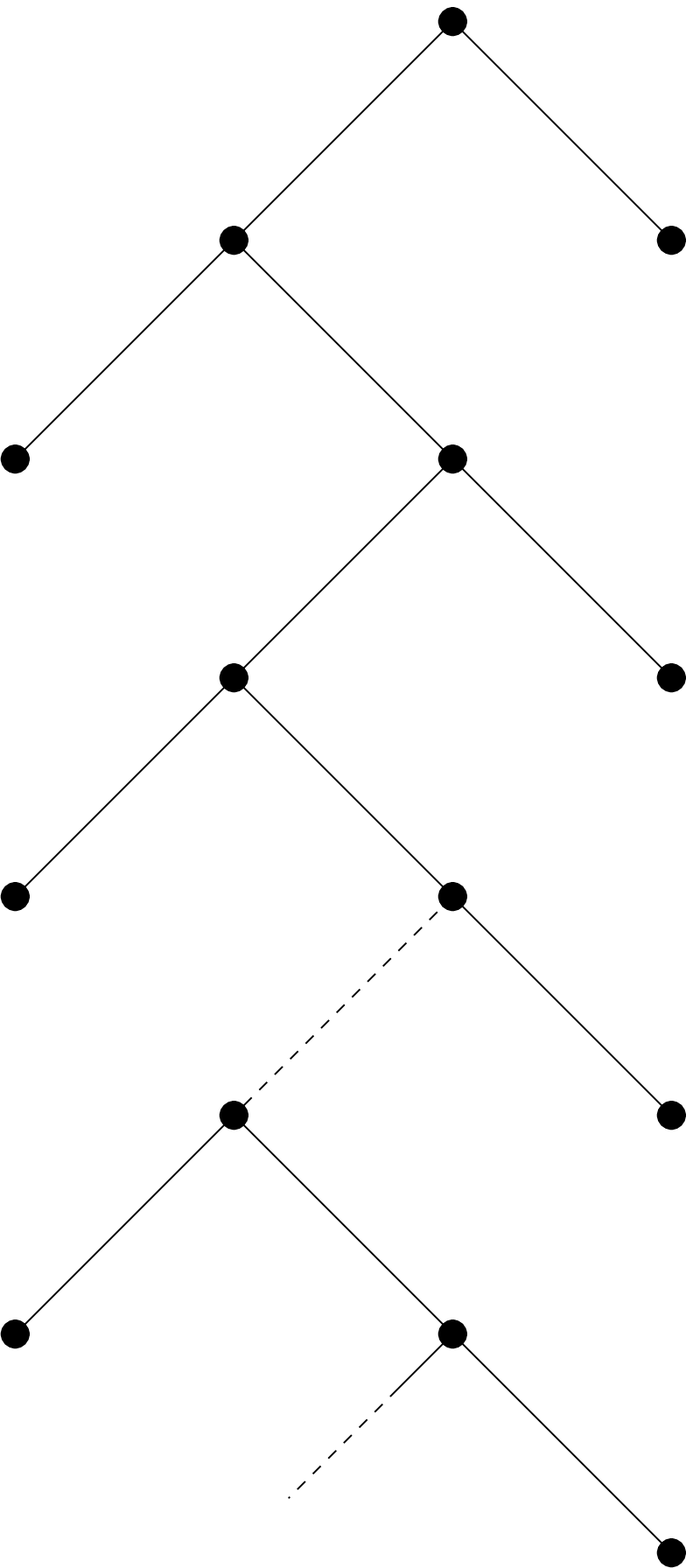}}
\put(220,10){\includegraphics[height=210pt]{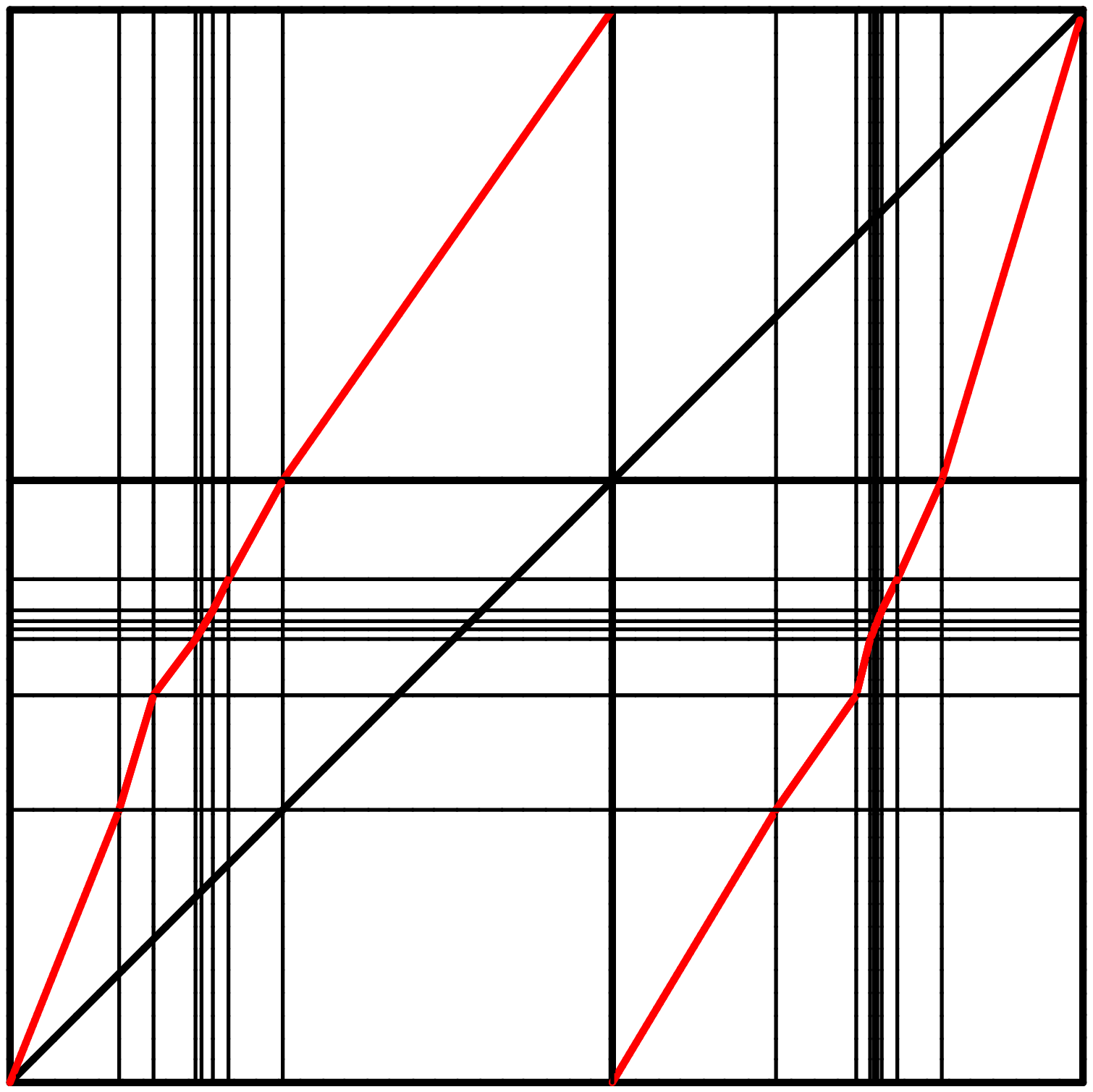}}

\thicklines
\put(116,158){$q_{1}$}
\put(32,133){$q_{00}$}
\put(112,106){$q_{011}$}
\put(28,80){$q_{0100}$}
\put(108,53){$q_{01011}$}
\put(24,27){$q_{(01)^n00}$}
\put(108,0){$q_{(01)^{n+1}1}$}
\put(63,10){$q_{(01)^\infty}$}

\put(244.4,13){\line(1,0){161.4}}
\put(284,20){$I_0$}
\put(364,20){$I_1$}
\put(244.4,10){\line(0,1){6}}
\put(335,10){\line(0,1){6}}
\put(405.8,10){\line(0,1){6}}

\put(220,34.2){\line(0,1){161.6}}
\put(225,77){$I_0$}
\put(225,154){$I_1$}
\put(217,34.2){\line(1,0){6}}
\put(217,125){\line(1,0){6}}
\put(217,195.8){\line(1,0){6}}
\end{picture}
\caption{\label{figBambou}
bamboo blossom probabilized context tree (on the left) and the associated dynamical system
(on the right).}
\end{figure}

The data of a corresponding VLMC consist in probability measures
on $\rond A$ indexed by the two families of finite contexts
  $$(q_{(01)^n1})_{n\geq 0}\mbox{\ and\ }(q_{(01)^n00})_{n\geq 0}$$
together with a probability measure on the infinite context $q_{(01)^{\infty}}$.

As before, assuming that $\pi$ is a stationary probability measure on $\rond L$, we
compute the probabilities of any $\pi(w)$, $\overline w$ being an internal node or
$\overline w$ being a context, as functions of the data and of both $\pi(1)$ and
$\pi(00)$.
Determination of stationary probabilities of cylinders based on both contexts $1$ and $00$  then leads to assumptions that guarantee existence and unicity of such a stationary probability measure.

\medskip

{\bf Computation of $\pi(w), \overline{w}$ context}

Two families of cylinders, namely $\rond L 1(10)^n$ and $\rond L 00(10)^n$, correspond to contexts.
For any $n\geq 0$,
$\pi(1(10)^{n+1})=\pi(1(10)^{n})q_{(01)^n1}(1)q_1(0)$
and
$\pi(00(10)^{n+1})=\pi(00(10)^{n})q_{(01)^n00}(1)q_1(0)$.
A straightforward induction implies thus that for any $n\geq 0$,
\begin{equation}
\label{probasContextesBambou}
\left\{
\begin{array}{l}
\pi(1(10)^{n})=\pi(1)c_n(1)
\\ \\
\pi(00(10)^{n})=\pi(00)c_n(00)
\end{array}
\right.
\end{equation}
where $c_0(1)=c_0(00)=1$ and
$$
\left\{
\begin{array}{l}
c_n(1)=\displaystyle q_1(0)^n\prod _{k=0}^{n-1}q_{(01)^k1}(1)
\\ \\
c_n(00)=\displaystyle q_1(0)^n\prod _{k=0}^{n-1}q_{(01)^k00}(1)
\end{array}
\right.
$$
for any $n\geq 1$.

\medskip

{\bf Computation of $\pi(w), \overline{w}$ internal node}

Two families of cylinders, $\rond L 0(10)^n$ and $\rond L (10)^n$, correspond to internal nodes. By disjoint union of events, they are related by
$$
\left\{
\begin{array}{l}
\pi(0(10)^{n})=\pi((10)^{n})-\pi(1(10)^{n})\\
\pi((10)^{n+1})=\pi(0(10)^{n})-\pi(00(10)^{n})
\end{array}
\right.
$$
for any $n\geq 0$.
By induction, this leads to: $\forall n\geq 0$,
\begin{equation}
\label{probasInternesBambou}
\left\{
\begin{array}{l}
\pi(0(10)^{n})=1-\pi(1)S_n(1)-\pi(00)S_{n-1}(00)
\\ \\
\pi((10)^{n})=1-\pi(1)S_{n-1}(1)-\pi(00)S_{n-1}(00)
\end{array}
\right.
\end{equation}
where $S_{-1}(1)=S_{-1}(00)=0$ and, for any $n\geq 0$,
$$
\left\{
\begin{array}{l}
S_n(1)=\sum _{k=0}^nc_k(1)
\\ \\
S_n(00)=\sum _{k=0}^nc_k(00).
\end{array}
\right.
$$
These formulae give, by quotients, the conditional probabilities on internal nodes defined by (\ref{probacondinternalnodes}) and appearing in
Formula~(\ref{probaStationnaireCylindres}).

\medskip

{\bf Computation of $\pi(1)$ and of $\pi(00)$}

\noindent
The context tree defines a partition of the set $\rond L$ of left-infinite sequences (see~(\ref{partition})).
In the case of bamboo blossom, this partition implies
\begin{eqnarray}
\label{equationPartitionBambou}
1-\pi((10)^\infty )
&=&\sum _{n\geq 0}\pi(1(10)^n)
+\sum _{n\geq 0}\pi(00(10)^n)\\
\label{equationPartitionBambou2}
&=& \sum_{n\geq 0}\pi(1)c_n(1) +  \sum_{n\geq 0}\pi(00)c_n(00).
\end{eqnarray}
We denote
$$
\left\{
\begin{array}{l}
S(1)=\sum _{n\geq 0}c_n(1)
\\ \\
S(00)=\sum _{n\geq 0}c_n(00) \in [1,+\infty].
\end{array}
\right.
$$
Note that the series $S(1)$ always converges. Indeed, the convergence is obvious if $q_1(0)\neq 1$; otherwise, $q_1(0) = 1$ and $q_1(1) = 0$, so that any $c_n(1)$, $n\geq 1$ vanishes and $S(1) = 1$. In the same way, the series $S(00)$ is finite as soon as $q_1(0)\neq 1$.

\begin{Prop}\label{mesureStationnaireBambou}
{\bf (Stationary measure on a bamboo blossom)}

Let $(U_n)_{n\geq 0}$ be a  VLMC defined by a probabilized bamboo blossom context tree.

\medskip

(i) Assume that $q_1(0)\neq 1$, then the Markov process $(U_n)_{n\geq 0}$ admits a stationary probability measure on $\rond L$ which is unique if and only if $S(1)-S(00)(1+q_1(0))\neq 0$.

\medskip

(ii) Assume that $q_1(0)=1$.

\smallskip

(ii.a) If $S(00)=\infty$, then $(U_n)_{n\geq 0}$ admits $\pi =\frac 12\delta _{(10)^\infty}+\frac 12\delta _{(10)^\infty 1}$ as unique stationary probability measure on $\rond L$.

\smallskip

(ii.b)  If $S(00)<\infty$, then $(U_n)_{n\geq 0}$ admits a one parameter family of stationary probability measures on $\rond L$.
\end{Prop}

\pff
{\it (i)} Assume that $q_1(0)\neq 1$ and that $\pi$ is a stationary probability measure.
Applying (\ref{probaStationnaireCalcul}) gives
\begin{equation}
\label{asymptotBambou}
\pi((10)^{\infty})=q_1(0)q_{(01)^{\infty}}(1)\pi((10)^{\infty})
\end{equation}
and consequently $\pi((10)^{\infty})=0$.
Therefore, Equation~(\ref{equationPartitionBambou}) becomes $S(1)\pi(1)+S(00)\pi(00)=1$.
We get another linear equation on $\pi(1)$ and $\pi(00)$ by disjoint
union of events:
$\pi(0)=1-\pi(1)
=\pi(10)+\pi(00)
=\pi(1)q_1(0)+\pi(00)$.
Thus $\pi(1)$ and $\pi(00)$ are solutions of the linear system
\begin{equation}
\label{systemeBambou}
\left\{
\begin{array}{l}
S(1)\pi(1)+S(00)\pi(00)=1
\\ \\
\left[ 1+q_1(0)\right]\pi(1)+\pi(00)=1.
\end{array}
\right.
\end{equation}
This system has a unique solution if and only if the determinantal assumption
$$
S(1)-S(00)\left[ 1+q_1(0)\right]\neq 0
$$
is fulfilled, which is a very light assumption (if this determinant happens to be zero, it suffices to modify one value of some $q_u$, $u$ context for the
assumption to be satisfied).
Otherwise, when the determinant vanishes, System~(\ref{systemeBambou}) is reduced to its second equation, so that it admits a one parameter family of solutions. Indeed, 
$$
1\leq S(1)\leq 1+q_1(0)(1-q_1(0))+\sum _{n\geq 2}q_1(0)^n(1-q_1(0))=1+q_1(0)
$$
and $S(00)\geq 1$, so that $S(1)-S(00)(1+q_1(0))=0$ implies that $S(1)=1+q_1(0)$ and $S(00)=1$.
In any case, System~(\ref{systemeBambou}) has at least one solution, which ensures the existence of a stationary probability measure with Formulae~(\ref{probasInternesBambou}), (\ref{probasContextesBambou}) and~(\ref{probaStationnaireCylindres}) by a standard argumentation. Assertions on unicity are straightforward.

\medskip

{\it (ii)} Assume that $q_1(0)=1$. This implies $q_1(1)=0$ and consequently $S(1)=1$. 
Thus, $\pi(1)$ and $\pi(00)$ are solutions of 
\begin{equation}
\label{systemeBambou2}
\left\{
\begin{array}{l}
\pi(1)+S(00)\pi(00)=1-\pi((10)^{\infty})
\\ \\
2 \pi(1)+\pi(00)=1.
\end{array}
\right.
\end{equation}
so that, since $S(00)\geq 1$, the determinantal condition $S(1)-S(00)(1+q_1(0))\neq 0$ is always fulfilled.

{\it (ii.a)} When $S(00)=\infty$, $\pi(00)=0$, $\pi(1)=\frac12$ and $\pi((10)^\infty )=\frac12$. This defines uniquely a stationary probability measure $\pi$.
Because of~(\ref{asymptotBambou}), $q_{(01)^{\infty}}(1)=1$ so that
$\pi ((10)^\infty 1)=\pi ((10)^\infty ))=\frac 12$.
This shows that $\pi =\frac 12\delta _{(10)^\infty}+\frac 12\delta _{(10)^\infty 1}$.

\medskip

{\it (ii.b)}  When $S(00)<\infty$, if we fix the value $a=\pi((10)^\infty )$, System (\ref{systemeBambou2}) has a unique solution that
determines in a unique way the stationary probability measure $\pi_a$.
\QED

\subsubsection{The associated dynamical system}
\label{sssec:blossomSD}

The vertical partition is made of the intervals $I_{(01)^n00}$ and $I_{(01)^n1}$ for $n\geq 0$. The horizontal partition
consists in the intervals $I_{0(01)^n00}$, $I_{1(01)^n00}$, $I_{0(01)^n1}$ and $I_{1(01)^n1}$ for $n\geq 0$, together with
the two intervals coming from the infinite context, namely $I_{0(01)^{\infty}}$ and $I_{1(01)^{\infty}}$.
If we make an hypothesis to ensure $\pi((10)^{\infty})=0$, then these two last intervals become two accumulation points
of the horizontal partition, $a_0$ and $a_1$. The respective positions of the intervals 
and the two accumulation points are given by the alphabetical order
  $$0(01)^{n-1}00<0(01)^n00<0(01)^{\infty}<0(01)^n1<0(01)^{n-1}1$$
  $$1(01)^{n-1}00<1(01)^n00<1(01)^{\infty}<1(01)^n1<1(01)^{n-1}1$$

\begin{Lem}
If $(q_{(01)^n00}(0))_{n\in{\N}}$ and $(q_{(01)^n1}(0))_{n\in{\N}}$ converge, then $T$ is right and left differentiable in $a_0$ and $a_1$ -- with possibly infinite derivatives -- and
  $$T'_\ell(a_0)=\lim_{n\to\infty}\frac{1}{q_{(01)^n00(0)}},\ \  T'_r(a_0)=\lim_{n\to\infty}\frac{1}{q_{(01)^n1(0)}}$$
  $$T'_\ell(a_1)=\lim_{n\to\infty}\frac{1}{q_{(01)^n00(1)}},\ \  T'_r(a_1)=\lim_{n\to\infty}\frac{1}{q_{(01)^n1(1)}}.$$
\end{Lem}

\pff
We use Lemma \ref{TAF}.
\QED

\subsubsection{Dirichlet series}

As for the infinite comb, the Dirichlet series of a source generated by a stationary bamboo blossom can be explicitly computed as a function of the SVLMC data.
For simplicity, we assume that the generic Condition {\it (i)} of Proposition~\ref{mesureStationnaireBambou}
is satisfied.
An internal node is of the form $(01)^n$ or $(01)^n0$ while a context writes $(01)^n00$ or $(01)^n1$.
Therefore, by disjoint union,
$$
\Lambda(s)=A(s)+\sum_{n\geq 0,w\in\rond W}\pi(w00(10)^n)^s
 +\sum_{n\geq 0,w\in\rond W}\pi(w1(10)^n)^s
 $$
 where
$$
A(s)=\sum_{n\geq 0}\pi((10)^n)^s+\sum_{n\geq 0}\pi(0(10)^n)^s
$$
is explicitly given by Formulae~(\ref{probasInternesBambou}) and~(\ref{systemeBambou}).
Because of the renewal property of the bamboo blossom, Formula~(\ref{probaStationnaireCalcul}) leads by
two straightforward inductions to $\pi (w00(10)^n)=\pi (w00)c_n(00)$ and $\pi (w1(10)^n)=\pi (w1)c_n(1)$
for any $n\geq 0$.
This implies that
$$
\Lambda (s)=A(s)+\Lambda_{00}(s)\sum_{n\geq 0}c_n^s(00)+\Lambda_{1}(s)\sum_{n\geq 0}c_n^s(1)
$$
where
$$
\Lambda _{00}(s)=\sum _{w\in\rond W}\pi(w00)^s
{\rm ~~and~~}
\Lambda _1(s)=\sum _{w\in\rond W}\pi(w1)^s.
$$
It remains to compute both Dirichlet series $\Lambda_{00}$ and $\Lambda_{1}$, which can be done by a
similar procedure.

\medskip
By disjoint union of finite words,
\begin{equation}
\label{lambda00}
\Lambda _{00}(s)=A_{00}(s)+\sum_{n\geq 0,w\in\rond W}\pi(w00(10)^n00)^s
+\sum_{n\geq 0,w\in\rond W}\pi(w1(10)^n00)^s
\end{equation}
where
$$
A_{00}(s)=\sum_{n\geq 0}\pi((10)^n00)^s+\sum_{n\geq 0}\pi(0(10)^n00)^s
$$
and
\begin{equation}
\label{lambda1}
\Lambda _1(s)=A_1(s)+\sum_{n\geq 0,w\in\rond W}\pi(w00(10)^n1)^s
+\sum_{n\geq 0,w\in\rond W}\pi(w1(10)^n1)^s
\end{equation}
with
$$
A_1(s)=\sum_{n\geq 0}\pi((10)^n1)^s+\sum_{n\geq 0}\pi(0(10)^n1)^s.
$$

\medskip
{\bf Computation of $A_1$ and $A_{00}$}

By disjoint union and Formula~(\ref{probaStationnaireCalcul}),
$$
\pi ((10)^{n+1}00)=\pi (0(10)^n00)-\pi (00(10)^n)q_{(01)^n00}(0)q_{00}(0),~n\geq 0
$$
and
$$
\pi (0(10)^n00)=\pi ((10)^n00)-\pi (1(10)^n)q_{(01)^n1}(0)q_{00}(0),~n\geq 1
$$
where $\pi (00(10)^n)$ and $\pi (1(10)^n)$ are already computed probabilities of contexts
(Formula~(\ref{probasContextesBambou})).
Since $\pi (000)=\pi (00)q_{00}(0)$, one gets recursively $\pi((10)^n00)$ and $\pi(0(10)^n00)$ from these two
relations as functions of the data.
This computes~$A_{00}$.
A very similar argument leads to an explicit form of $A_1$.

\vfill\eject
{\bf Ultimate computation of $\Lambda _1$ and $\Lambda _{00}$}

Start with~(\ref{lambda00}) and~(\ref{lambda1}).
 As above,  for any $n\geq 0$, by induction and with Formula~(\ref{probaStationnaireCalcul}),
 $$
 \pi (w00(10)^n00)=\pi (w00)c_n(00)q_{(01)^n00}(0)q_{00}(0).
 $$
 In the same way, but only when $n\geq 1$,
 $$
 \pi (w1(10)^n00)=\pi (w1)c_n(1)q_{(01)^n1}(0)q_{00}(0).
 $$
 Similar computations lead to similar formulae for 
 $\pi (w00(10)^n1)$ and $\pi (w1(10)^n1)$, for any $n\geq 0$.
 So,~(\ref{lambda00}) and~(\ref{lambda1}) lead to
 \begin{equation}
 \label{lambda00Split}
 \Lambda _{00}(s)=A_{00}(s)+\Lambda _{100}(s)
 +\Lambda _{00}(s)B_{00}(s)
 +\Lambda _{1}(s)B_1(s)
 \end{equation}
 where $B_{00}(s)$ and $B_1(s)$ are explicit functions of the data and where
 $$
 \Lambda _{100}(s)=\sum _{w\in\rond W}\pi (w100).
 $$
 As above, after disjoint union of words, splitting by Formula~(\ref{probaStationnaireCalcul}) and double
 induction, one gets
 $$
 \Lambda _{100}(s)=A_{100}(s)+\Lambda _{00}(s)C_{00}(s)+\Lambda _{1}(s)C_1(s)
 $$
 where $A_{100}(s)$, $C_{00}(s)$ and $C_{1}(s)$ are explicit series, functions of the data.
 Replacing $\Lambda _{100}$ by this value in Formula~(\ref{lambda00Split}) leads to a first linear equation
 between $\Lambda _{00}(s)$ and $\Lambda _1(s)$.
 A second linear equation between them is obtained from~(\ref{lambda1}) by similar arguments.
 Solving the system one gets with both linear equations gives an explicit form of $\Lambda _{00}(s)$ and
 $\Lambda _1(s)$ as functions of the data, completing the expected computation.

\subsubsection{Generating function for the exact distribution of word occurrences in a sequence generated by a bamboo blossom}\label{ssec:generatrice-bambou}

Let us consider the process $X=(X_n)_{n\geq 0}$ of final letters of $(U_n)_{n\geq 0}$ in the particular case of a SVLMC defined by a bamboo blossom. We only deal with finite words $w$ such that $\overline w$ is not an internal node, \ie $\overline w$ is a finite context or $\overline w\notin \rond T$. One can see that such a word of length $k>1$ can be written in the form $*11(10)^{\ell}1^{p}$ or $*00(10)^{\ell}1^{p}$, with $p\in \{0,1\}$ and $\ell \in \N$, where $*$ stands for any finite word. 
\begin{Prop}
\label{propgenbamboo} 
For a SVLMC defined by a bamboo blossom, with notations of Section~\ref{ssec:generatrice-peigne}, the generating function of the first occurrence of a finite word $w=w_1 \ldots w_k$ is given for $|x|<1$ by
$$
\Phi_w^{(1)}(x)=\frac{x^k\pi(w)}{(1-x)S_w(x)}
$$
and the generating function of the $r$\textsuperscript{th} occurrence of $w$ is given by
$$
\Phi_w^{(r)}(x)=\Phi_w^{(1)}(x)\left(1-\frac{1}{S_w(x)}\right)^{r-1},
$$
where

{\bf (i)} if $w$ is of the form $*00(10)^{\ell}$ or $*11(01)^{\ell}0$, with $\ell \in\g N$,
$S_w(x)$ is defined in Proposition~\ref{propgengen} and

{\bf (ii)} if $w$ is of the form $*00(10)^{\ell}1$, $\ell \in\g N$,
\begin{eqnarray*}
S_w(x)&=& C_w(x) + \sum_{j=k}^{\infty}q_{1(01)^{\ell}00}^{(j)}(w)x^j,\\
C_w(x)&=& 1+\sum_{j=1}^{k-1}\ind{w_{j+1}\ldots w_{k}=w_{1}\ldots w_{k-j}}q_{1(01)^{\ell}00}^{(j)}\left(w_{j+1} \ldots w_k\right)x^j.
\end{eqnarray*}
and if $w$ is of the form $*11(01)^{\ell}$, $\ell \in\g N$,
\begin{eqnarray*}
S_w(x)&=& C_w(x) + \sum_{j=k}^{\infty}q_{(10)^{\ell}11}^{(j)}(w)x^j,\\
C_w(x)&=& 1+\sum_{j=1}^{k-1}\ind{w_{j+1}\ldots w_{k}=w_{1}\ldots w_{k-j}}q_{(10)^{\ell}11}^{(j)}\left(w_{j+1} \ldots w_k\right)x^j.
\end{eqnarray*}
\end{Prop}

\pff
{\bf (i)} We first deal with the words $w$ such that 
\[\lpref(w)=(01)^{\ell}00 \quad \mbox{or} \quad \lpref(w)=(01)^{\ell}1.\] 
Let us denote by $p_n$ the probability that $T_w^{(1)}=n$. Proceeding exactly in the same way as for Proposition~\ref{propgengen}, from the decomposition 
\[\pi(w)=p_n+\sum_{z=k}^{n-1}p_z\Proba\left(X_{n-k+1}\ldots X_{n}=w|T_w^{(1)}=z\right),\]
and due to the renewal property of the bamboo, one has
\begin{eqnarray*}
\pi(w)=p_n&+&\sum_{z=k}^{n-k}p_z\Proba\left(U_n\in \rond L w\Big|U_z \in \rond L\suff(w)\right)\\
&+&\sum_{z=n-k+1}^{n-1}p_z\ind{w_{n-z+1}\ldots w_k=w_1 \ldots w_{k-n+z}}\Proba\left(U_n\in \rond L w\Big|U_z \in \rond L\suff(w)\right)
\end{eqnarray*}
where $\suff (w)$ is the suffix of $w$ equal to the reversed word of $\lpref (w)$.
Hence, for $x<1$, it comes
\begin{eqnarray*}
\frac{x^k\pi(w)}{1-x}=\sum_{n=k}^{+\infty}p_nx^n&+&\sum_{n=k}^{+\infty}x^n\sum_{z=k}^{n-k}p_zq_{\petitlpref(w)}^{(n-z)}(w)\\
&+& \sum_{n=k}^{+\infty}x^n\sum_{z=n-k+1}^{n-1}p_z\ind{w_{n-z+1}\ldots w_k=w_1 \ldots w_{k-n+z}}q_{\petitlpref(w)}^{(n-z)}(w)
\end{eqnarray*}
which leads to the expression of $\Phi_w^{(1)}(x)$ given in Proposition~\ref{propgengen}. The $r$\textsuperscript{th} occurrence can be derived exactly in the same way from the decomposition
\[\{w \mbox{ at }n\}=\{T_w^{(1)}=n\}\cup\{T_w^{(2)}=n\}\cup\ldots\cup\{T_w^{(r)}=n\}\cup\{T_w^{(r)}<n \mbox{ and }w \mbox{ at }n\}.\]

{\bf (ii)} In the particular case of words $w=*00(10)^{\ell}1$, the main difference is that the context $1$ is not sufficient for the renewal property. The computation relies on the equality
\[\Proba\left(X_{n-k+1}\ldots X_{n}=w|T_{w}^{(1)}=z\right)=\Proba\left(X_{n-k+1}\ldots X_n=w|X_{z-2\ell-2}\ldots X_z=00(10)^{\ell}1\right).\]
The sketch of the proof remains the same replacing $q_{\petitlpref(w)}(w)$ by $q_{1(01)^{\ell}00}(w)$. The case $w=*11(01)^{\ell}$ is analogous.
\QED

\section{Some remarks, extensions and open problems} 
\label{ouverture}
\subsection{Stationary measure for a general VLMC}
\label{secStationaryMeasureGeneralVlmc}

Infinite comb and bamboo blossom are two instructive but very particular examples, close to renewal processes. Nevertheless, we think that an analogous of Proposition \ref{mesureStationnairePeigne} or Proposition \ref{mesureStationnaireBambou} can be written for a VLMC defined by a general context tree with a finite or denumerable number of infinite branches. 

In order to generalize the proofs, it is clear that Formula (\ref{probaStationnaireCylindres}) in Lemma~\ref{lemProbaStationnaireCylindres} is crucial. In this formula, for a given finite word $w=\alpha _1\dots\alpha _N\in\rond W$ it is important to check whether the subwords $\lpref (\alpha _1\dots \alpha _k), k <N,$ are   internal nodes of the tree or not. Consequently, the following concept of \emph{minimal context} is natural.

\begin{Def} {\bf (Minimal context)}
Define the following binary relation on the set of the finite contexts as follows: 
$$\forall u,v \in \rond C^{F}, \hskip 5pt  u\prec v\Longleftrightarrow\exists w,w'\in\rond W,~v=wuw'$$
(in other words $u$ is a sub-word of $v$). This relation is a partial order. In  a context tree, a finite context is called \emph{minimal} when it is minimal for this partial order on contexts. 
\end{Def}

\begin{Rem} {\bf (Alternative definition of a minimal context)}
Let $\rond T$ be a context tree.
Let $c =  \alpha_N\ldots \alpha_1$ be a finite context of $\rond T$.
Then $c$ is minimal if and only if
$\forall k\in \{1, \ldots, N-1\}, \lpref(\alpha_1 \ldots \alpha_k) \notin \rond C^F(\rond T)$.
\end{Rem}

\begin{Ex}
In the infinite comb, the only minimal context is $1$. In the bamboo blossom, the minimal contexts are $1$ and $00$.
\end{Ex}

\begin{Rem}

There exist some trees with infinitely many infinite leaves and a finite number of minimal contexts. Take the infinite comb and at each $0^k$ branch another infinite comb. In such a tree, the finite leaf\ $10$ is the only minimal context.

Nonetheless, a tree with a finite number of infinite contexts has necessarily a finite number of minimal contexts.
\end{Rem}

As one can see for the infinite comb or for the bamboo blossom (see Sections~\ref{sssec:combSVLMC}
and~\ref{sssec:blossomSVLMC}), minimal contexts play a special role in the computation of stationary
probability measures.
First of all, when $\pi$ is a stationary probability measure and $w$ a finite word such that $\overline w\notin\rond T$, Formula~(\ref{probaStationnaireCylindres})
implies that $\pi (w)$ is a rational monomial of the data $q_c(\alpha )$ and of the $\pi (u)$ where $\overline u$
belongs to $\rond T$.
This shows that any stationary probability is determined by its values on the nodes of the context tree.
In both examples, we compute these values as functions of the data and
of the $\pi (m)$, where $\overline m$ are minimal contexts, and we finally write a  rectangular linear system
satisfied by these $\pi (m)$.
Assuming that this system has maximal rank can be viewed as making an irreducibility condition for the Markov chain on $\rond L$.
We conjecture that this situation happens in any case of VLMC.

In the following example, we detail the above procedure, in order to understand how the two main principles (the partition (\ref{partition}) and the disjoint union) give the linear system leading to the irreducibility condition.

\begin{Ex}
Let $\rond T$ be a probabilized context tree corresponding to Figure~\ref{figPeigneFini}
(finite comb with $n+1$ teeth).
There are two minimal contexts: $1$ and $0^{n+1}$.
\begin{figure}[!h]
\begin{picture}(400,150)
\put(145,15){\includegraphics[width=150pt]{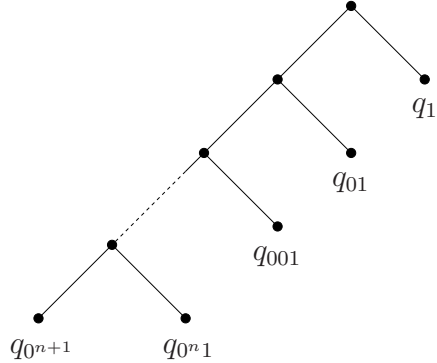}}
\put(136,4){$q_{0^{n+1}}$}
\put(194,4){$q_{0^n1}$}
\put(228,39){$q_{001}$}
\put(258,67){$q_{01}$}
\put(288,95){$q_{1}$}
\end{picture}
\caption{\label{figPeigneFini}$(n+1)$-teeth comb probabilized context tree.}
\end{figure}
Assume that $\pi$ is a stationary probability measure on $\rond L$.
Like in the case of the infinite comb, the probability of a word that corresponds to a teeth is 
$\pi (10^k)=\pi (1)c_k$, $0\leq k\leq n$ where $c_k$ is the product defined by~(\ref{definitioncnPeigne}).
Also, the probabilities of the internal nodes and of the handle are 
$$
\pi (0^k)=1-\pi (1)S_{k-1},~~ 0\leq k\leq n+1,
$$
where $S_p:=\sum _{j=0}^{p}c_j$.
By means of these formulae, $\pi$ is determined by $\pi (1)$.

In order to compute $\pi (1)$, one can proceed as follows.
First, by the partition principle (\ref{partition}), we have $1=\pi (0^{n+1})+\pi (1)\sum _{k=0}^{n}c_k$.
Secondly, by disjoint union,
$$\pi (0^{n+1})=\pi (0^{n+2})+\pi (10^{n+1})=\pi (0^{n+1})q_{0^{n+1}}(0)+\pi (10^n)q_{0^n1}(0).$$
This implies the linear relation between both minimal contexts probabilities:
$$
\left\{
\begin{array}{ll}
\pi (0^{n+1})+S_n\pi (1)=1
\\
q_{0^{n+1}}(1)\pi (0^{n+1})-q_{0^n1}(0)c_n\pi (1)=0.
\end{array}
\right.
$$
In particular, this leads to the irreducibility condition $q_{0^{n+1}}(1)S_n+q_{0^n1}(0)c_n\neq 0$
for the VLCM to admit a stationary probability measure.
One can check that this irreducibility condition is the classical one for the corresponding $\rond A$-valued
Markov chain of order $n+1$.
\end{Ex}
\begin{Ex}
\label{bb4f}
Let $\rond T$ be a probabilized context tree corresponding to Figure~\ref{figArbreMapBambou4Fleurs} (four flower bamboo). This tree provides another example of  computation procedure using Formulae~(\ref{probaStationnaireCalcul}) and~(\ref{probaStationnaireCylindres}), the partition principle (\ref{partition}) and the disjoint union. This VLMC admits a unique stationary probability measure if the determinantal condition 
\[q_{00}(1)[1+q_1(0)]+q_1(0)^2q_{010}(0)+q_1(0)q_1(1)q_{011}(0)\neq 0\]
is satisfied; it is fulfilled if none of the $q_c$ is trivial. 
\end{Ex}
\subsection{Tries}

In a first kind of problems, $n$  words independently produced by a source are inserted in a trie. There are results on the classical parameters of the trie (size, height, path length) for a dynamical source (\cite{ClementFlajoletVallee}), which rely on the existence of a spectral gap for the underlying dynamical system. We would like to extend these results to cases when there is no spectral gap, as may be guessed in the infinite comb example.

\medskip

Another interesting application consists in producing a suffix trie from \emph{one} sequence coming from a VLMC dynamical source, and analyzing its parameters. 
For his analysis, \cite{Szp} puts some mixing assumptions (called strong $\alpha$-mixing) on the source.
A first direction consists in trying to find the mixing type of a VLMC dynamical source. In a  second direction, we plan to use the generating function for the occurrence of words to improve these results.

\subsection*{Acknowledgements}

We are very grateful to Antonio Galves, who introduced us to the challenging VLMC topics. We warmly thank Brigitte Vall\'ee for valuable and  stormy discussions.

\bibliographystyle{plainnat} 
\bibliography{vomc-sysdyn}

\eject
{\fontauthors
Peggy C\'enac}

Universit\'e de Bourgogne

Institut de Math\'ematiques de Bourgogne

IMB UMR 5584 CNRS

9 rue Alain Savary - BP 47870, 21078 DIJON CEDEX

{\fontemail peggy.cenac@u-bourgogne.fr}

\vskip 15pt
{\fontauthors
Brigitte Chauvin}

INRIA Rocquencourt, project Algorithms,

Domaine de Voluceau B.P.105, 78153 Le Chesnay CEDEX (France),

and

Laboratoire de Math\'ematiques de Versailles,

CNRS, UMR 8100,

Universit\'e de Versailles - St-Quentin,

45, avenue des Etats-Unis, 78035 Versailles CEDEX (France)

{\fontemail chauvin@math.uvsq.fr}

\vskip 15pt
{\fontauthors
F\'ed\'eric Paccaut}

LAMFA

CNRS, UMR 6140

Universit\'e  de Picardie Jules Verne

33, rue Saint-Leu, 80039 Amiens, France

{\fontemail frederic.paccaut@u-picardie.fr}

\vskip 15pt
{\fontauthors
Nicolas Pouyanne}

Laboratoire de Math\'ematiques de Versailles,

CNRS, UMR 8100,

Universit\'e de Versailles - St-Quentin,

45, avenue des Etats-Unis, 78035 Versailles CEDEX (France)

{\fontemail pouyanne@math.uvsq.fr}

\end{document}